\documentclass[12pt]{article}
\usepackage{amsfonts,amsmath}
\usepackage{epic}

\renewcommand{\appendix}{%
\renewcommand{\section}{%
\newpage
\thispagestyle{plain}%
\secdef\Appendix\sAppendix}%
\setcounter{section}{0}%
\renewcommand{\thesection}{\Alph{section}}%
}
\newcommand{\Appendix}[2][?]{%
\refstepcounter{section}%
\addcontentsline{toc}{appendix}%
{\protect\numberline{\appendixname~\thesection}#1}%
{\flushleft\LARGE\bfseries\appendixname\ \thesection\par
\centering#2\par}%
\sectionmark{#1}\vspace{\baselineskip}}

\newcommand{\sAppendix}[1]{%
{\flushright\large\bfseries\appendixname\par
\centering#1\par}%
\vspace{\baselineskip}}

\def\be{\begin{equation}}

\def\bea{\begin{eqnarray}}

\def\eea{\end{eqnarray}}

\begin{document}

\pagestyle{empty}

\begin{center}

\textsf{\Huge {\bf A new class ${\hat o}_N$ of statistical models:
Transfer matrix eigenstates, chain Hamiltonians, factorizable
$S$-matrix}}

\vspace{10mm}

{\large B. Abdesselam$^{a,}$\footnote{Email:
boucif@cpht.polytechnique.fr and  boucif@yahoo.fr} and A.
Chakrabarti$^{b,}$\footnote{Email: chakra@cpht.polytechnique.fr}}

  \vspace{5mm}

  \emph{$^a$ Laboratoire de Physique Quantique de la
Mati\`ere et de Mod\'elisations Math\'ematiques, Centre
Universitaire de Mascara, 29000-Mascara, Alg\'erie\\
and \\
Laboratoire de Physique Th\'eorique, Universit\'e d'Oran
Es-S\'enia, 31100-Oran, Alg\'erie}
  \\
  \vspace{3mm}
  \emph{$^b$ Centre de Physique Th{\'e}orique, CNRS UMR 7644}
  \\
  \emph{Ecole Polytechnique, 91128 Palaiseau Cedex, France.}
  \vspace{3mm}

\end{center}

\begin{abstract}
{\small \noindent Statistical models corresponding to a new class
of braid matrices ($\hat{o}_N;\,N\geq 3$) presented in a previous
paper are studied. Indices labeling states spanning the $N^r$
dimensional base space of $T^{(r)}\left(\theta\right)$, the $r$-th
order transfer matrix are so chosen that the operators $W$ (the
sum of the state labels) and (CP) (the circular permutation of
state labels) commute with $T^{(r)}\left(\theta\right)$. This
drastically simplifies the construction of eigenstates, reducing
it to solutions of relatively small number of simultaneous linear
equations. Roots of unity play a crucial role. Thus for
diagonalizing the 81 dimensional space for $N=3$, $r=4$, one has
to solve a maximal set of 5 linear equations. A supplementary
symmetry relates invariant subspaces pairwise ($W=(r,Nr)$ and so
on) so that only one of each pair needs study. The case $N=3$ is
studied fully for $r=\left(1,2,3,4\right)$. Basic aspects for all
$\left(N,r\right)$ are discussed. Full exploitation of such
symmetries lead to a formalism quite different from, possibly
generalized, algebraic Bethe ansatz. Chain Hamiltonians are
studied. The specific types of spin flips they induce and
propagate are pointed out. The inverse Cayley transform of the YB
matrix giving the potential leading to factorizable $S$-matrix is
constructed explicitly for $N=3$ as also the full set of
$\hat{R}tt$ relations. Perspectives are discussed in a final
section.}
\end{abstract}

\rightline{math.QA/0607379}

\pagestyle{plain} \setcounter{page}{1}

\section{Introduction}
\setcounter{equation}{0}

New classes of braided matrices were presented in recent papers
\cite{R1,R2}. Statistical models corresponding to \cite{R1} have
been presented in \cite{R3}. Here we present those corresponding
to \cite{R2}. Different types of statistical models thus obtained
will be compared at the end (sec. 7). In \cite{R2} two distinct
classes of braid matrices $\left({\hat o}_N,{\hat p}_N\right)$
were presented. Here we consider only the ${\hat o}_N$
$\left(N\geq 3\right)$. For real, positive values of the parameter
$q$ and a certain domain (depending on $q$ and $N$) of the
spectral parameter $\theta$, one obtains $N^2\times N^2$ braid
matrices with all nonzero elements real, positive giving
nonnegative Boltzmann weights. For the class ${\hat p}_N$ one
encounters both positive and negative elements and thus one would
need suitable reinterpretation of the corresponding Boltzmann
weights.

We first recapitulate briefly the ${\hat o}_N$ braid matrices
\cite{R2}. The $N^2\times N^2$ baxterized braid matrices
satisfying (in standard notations)
\begin{equation}
{\hat R}_{12}\left(\theta\right){\hat
R}_{23}\left(\theta+\theta'\right){\hat
R}_{12}\left(\theta'\right)={\hat R}_{23}\left(\theta'\right){\hat
R}_{12}\left(\theta+\theta'\right){\hat R}_{23}\left(\theta\right)
\end{equation}
are given by
\begin{equation}
{\hat R}\left(\theta\right)=I-\frac{\sinh
\theta}{\sinh\left(\eta+\theta\right)}P_0',
\end{equation}
where
\begin{equation}
e^\eta+e^{-\eta}=\left[N-1\right]+1\equiv
\frac{q^{N-1}-q^{-N+1}}{q-q^{-1}}+1
\end{equation} and
\begin{equation}
P_0'=\sum_{i,j=1}^Nq^{\rho_{j'}-\rho_j}\left(ij\right)\otimes
\left(i'j'\right)
\end{equation}
with the follwing notations:
\begin{enumerate}
\item The $N\times N$ matrix $\left(ij\right)$ has only one
non-zero element, unity, on row $i$ and column $j$ and
\begin{equation}
\left(i',j'\right)=\left(N-i+1,N-j+1\right).
\end{equation}

\item The $N$-tuple $\left(\rho_1,\rho_2,\ldots,\rho_N\right)$ is defined as
\begin{equation}
\left(n-\frac 12,n-\frac 32,\ldots,\frac 12,0,-\frac 12,\ldots,
-n+\frac 12\right)
\end{equation}
for $N=2n+1$ and
\begin{equation}
\left(n-1,n-2,\ldots,1,0,0,-1,\ldots, -n+1\right)
\end{equation}
for $N=2n$.
\end{enumerate}

Of the three projectors $\left(P_+,P_-,P_0\right)$ providing a
spectral resolution of $SO_q\left(N\right)$ braid matrices only
\begin{equation}
P_0'=\left(\left[N-1\right]+1\right)P_0
\end{equation}
appears in our class. To signal this provenance (along with
crucial differences) our class is designated as ${\hat o}_N$. More
relevant discussions can be found in \cite{R2}.

We now introduce the permutation matrix
\begin{equation}
P=\sum_{i,j}\left(ij\right)\otimes \left(ji\right),\qquad \qquad
P^2=I
\end{equation}
the Yang-Baxter matrix
\begin{equation}
R\left(\theta\right)=P{\hat R}\left(\theta\right)
\end{equation}
and the monodromy matrices satisfying
\begin{equation}
{\hat
R}\left(\theta-\theta'\right)\left(t\left(\theta\right)\otimes
t\left(\theta'\right)\right)=\left(t\left(\theta'\right)\otimes
t\left(\theta\right)\right){\hat R}\left(\theta-\theta'\right)
\end{equation}
The $t$-matrix satisfying (1.11) is $N\times N$ in terms of the
blocks
\begin{equation}
t_{ij},\qquad \left(i,j=1,\ldots,N\right)
\end{equation}
each $t_{ij}$ being itself a matrix whose dimension is prescribed
as follows. One starts with $N\times N$ blocks $t_{ij}$ obtained
from the standard prescription (satisfying (1.11))
\begin{equation}
t^{(1)}\left(\theta\right)=P\hat{R}\left(\theta\right)=R\left(\theta\right)
\end{equation}
and then a hierarchy is obtained implementing the coproduct
prescription
\begin{equation}
t^{(r)}_{ij}\left(\theta\right)=\sum_{k_1,\ldots,k_{r-1}}t^{(1)}_{ik_1}\left(\theta\right)\otimes
t^{(1)}_{k_1k_2}\left(\theta\right)\otimes\cdots \otimes
t^{(1)}_{k_{r-1},j}\left(\theta\right).
\end{equation}
Starting with (1.13), this prescription assures that
$t^{(r)}\left(\theta\right)$ satisfies (1.11).

Now the transfer matrix is defined, for each order $r$, as
\begin{equation}
T^{(r)}\left(\theta\right)=\sum_{i=1}^N
t^{(r)}_{ii}\left(\theta\right).
\end{equation}
The trace and more generally the eigenstates and the eigenvalues
of $T^{(r)}\left(\theta\right)$ provide crucial properties of the
statistical mechanical model associated with
$\hat{R}\left(\theta\right)$. In particular, (1.1), (1.11),
(1.13), (1.14), (1.15) all together assure the commutativity
\begin{equation}
\left[T\left(\theta\right),T\left(\theta'\right)\right]=0.
\end{equation}
Commutative transfer matrices provide the crucial feature of
exactly solvable models of statistical mechanics, the braid
matrices encoding star-triangle relations \cite{R4}. For our
specific case $\left(\hat{o}_N\right)$ we illustrate, in the
following section, some basic features for the simplest case
($N=3$). Certain aspects for $N>3$ will be presented afterwards
(sec. 5).

Define
\begin{equation}
K\left(\theta\right)=-\frac{\sinh\theta}{\sinh\left(\eta+\theta\right)},
\end{equation}
where (setting $N=3$ in (1.3))
\begin{equation}
e^\eta+e^{-\eta}=q+q^{-1}+1.
\end{equation}
For
\begin{equation}
-\eta<\theta<0,\qquad K\left(\theta\right)>0.
\end{equation}
For
\begin{eqnarray}
&&\theta=0,\qquad K\left(0\right)=0,\nonumber\\
&& \theta=-\frac\eta 2,\qquad K\left(-\frac\eta 2\right)=1,\nonumber\\
&& \theta\longrightarrow -\eta ,\qquad
K\left(\theta\right)\longrightarrow +\infty.
\end{eqnarray}
Henceforward we consider the domain (1.19).

\section{Trace of the transfer matrix from iterative structure}
\setcounter{equation}{0}

 The standard prescription (1.13) yields
for $\hat{o}_3$
\begin{equation}
t^{(1)}\left(\theta\right)=P\hat{R}\left(\theta\right)=P\left(I+K\left(\theta\right)P_0'\right)
\end{equation}
and hence (suppressing now the argument $\theta$ for simplicity)
\begin{eqnarray}
&&t_{11}^{(1)}=\left|\begin{array}{ccc}
   1 & 0 & 0 \\
   0 & 0 & 0 \\
   0 & 0 & K \\
\end{array}\right|,\qquad t_{12}^{(1)}=\left|\begin{array}{ccc}
   0 & 0 & 0 \\
   1 & 0 & 0 \\
   0 & q^{1/2}K & 0\\
\end{array}\right|,\qquad
t_{13}^{(1)}=\left|\begin{array}{ccc}
   0 & 0 & 0 \\
   0 & 0 & 0 \\
   1+qK & 0 & 0 \\
\end{array}\right|,\nonumber\\
&&t_{21}^{(1)}=\left|\begin{array}{ccc}
   0 & 1 & 0 \\
   0 & 0 & q^{-1/2}K \\
   0 & 0 & 0 \\
\end{array}\right|,\qquad t_{22}^{(1)}=\left|\begin{array}{ccc}
   0 & 0 & 0 \\
   0 & 1+K & 0 \\
   0 & 0 & 0 \\
\end{array}\right|,\qquad
t_{23}^{(1)}=\left|\begin{array}{ccc}
   0 & 0 & 0 \\
   q^{1/2}K & 0 & 0 \\
   0 & 1 & 0 \\
\end{array}\right|,\\
&&t_{31}^{(1)}=\left|\begin{array}{ccc}
   0 & 0 & 1+q^{-1}K \\
   0 & 0 & 0 \\
   0 & 0 & 0 \\
\end{array}\right|,\qquad
t_{32}^{(1)}=\left|\begin{array}{ccc}
   0 & q^{-1/2}K & 0 \\
   0 & 0 & 1 \\
   0 & 0 & 0 \\
\end{array}\right|,\qquad
t_{33}^{(1)}=\left|\begin{array}{ccc}
   K & 0 & 0 \\
   0 & 0 & 0 \\
   0 & 0 & 1 \\
\end{array}\right|.\nonumber
\end{eqnarray}
All $\theta$-dependence is contained in the parameter $K$, as
defined in (1.17). Starting with the $3\times 3$ blocks the
prescription (1.14) gives $3^r\times 3^r$ blocks $t_{ij}^{(r)}$.
The recursion relations for our case ($N=3$) are (for $j=1,2,3$)
\begin{eqnarray}
&&t_{1j}^{(r+1)}=\left|\begin{array}{ccc}
   t_{1j}^{(r)} & 0 & 0 \\
   t_{2j}^{(r)} & 0 & 0 \\
   \left(1+qK\right)t_{3j}^{(r)} & q^{1/2}Kt_{2j}^{(r)} & Kt_{1j}^{(r)} \\
\end{array}\right|,\\
&&t_{2j}^{(r+1)}=\left|\begin{array}{ccc}
   0 & t_{1j}^{(r)} & 0 \\
   q^{1/2}K t_{3j}^{(r)} & \left(1+K\right)t_{2j}^{(r)} & q^{-1/2}Kt_{1j}^{(r)} \\
   0 & t_{3j}^{(r)} & 0 \\
\end{array}\right|,\\
&&t_{3j}^{(r+1)}=\left|\begin{array}{ccc}
   Kt_{3j}^{(r)} & q^{-1/2}Kt_{2j}^{(r)} & \left(1+q^{-1}K\right)t_{1j}^{(r)} \\
   0 & 0 & t_{2j}^{(r)} \\
   0 & 0 & t_{3j}^{(r)}  \\
\end{array}\right|.
\end{eqnarray}
The transfer matrix is iterated as
\begin{eqnarray}
&&T^{(r+1)}=t_{11}^{(r+1)}+t_{22}^{(r+1)}+t_{33}^{(r+1)}\nonumber\\
&&\phantom{T^{(r+1)}}=\left|\begin{array}{ccc}
   t_{11}^{(r)}+Kt_{33}^{(r)} & t_{12}^{(r)}+q^{-1/2}Kt_{23}^{(r)} & \left(1+q^{-1}K\right)t_{13}^{(r)} \\
   t_{21}^{(r)}+q^{1/2}Kt_{32}^{(r)} & \left(1+K\right)t_{22}^{(r)} & q^{-1/2}Kt_{12}^{(r)}+t_{23}^{(r)} \\
   \left(1+qK\right)t_{31}^{(r)} & q^{1/2}Kt_{21}^{(r)}+t_{32}^{(r)} &Kt_{11}^{(r)}+t_{33}^{(r)} \\
\end{array}\right|.
\end{eqnarray}
Hence
\begin{eqnarray}
&&\hbox{Tr}\left(T^{(r+1)}\right)=\hbox{Tr}\left(t_{11}^{(r)}+Kt_{33}^{(r)}+\left(1+K\right)t_{22}^{(r)}+
Kt_{11}^{(r)}+t_{33}^{(r)}\right)\nonumber\\
&&\phantom{\hbox{Tr}\left(T^{(r+1)}\right)}=\left(1+K\right)\hbox{Tr}\left(t_{11}^{(r)}+t_{22}^{(r)}+t_{33}^{(r)}\right)
\nonumber\\
&&\hbox{Tr}\left(T^{(r+1)}\right)=
\left(1+K\right)\hbox{Tr}\left(T^{(r)}\right).
\end{eqnarray}
But from (2.2),
\begin{equation}
\hbox{Tr}\left(T^{(1)}\right)=\hbox{Tr}\left(t_{11}^{(1)}+t_{22}^{(1)}+t_{33}^{(1)}\right)=3\left(1+K\right).
\end{equation}
Hence
\begin{equation}
\hbox{Tr}\left(T^{(r)}\right)=3\left(1+K\right)^r.
\end{equation}
Thus we obtain the trace of $T^{(r)}$ for all $r$ directly without
constructing explicitly the eigenstates and the $3^r$ eigenvalues.
But the latter being of crucial interest we now turn to their
systematic explicit constructions.

\section{Eigenstates and eigenvalues $(N=3)$}
\setcounter{equation}{0}

For $N=3$ the transfer matrix $T^{(r)}\left(\theta\right)$ of
order $r$ acts on a space of dimension $3^r$. Construction of
eigenstates corresponds to diagonalization of $T^{(r)}$ on such a
base space. But basic symmetries of $T^{(r)}$ (Sec. 1) for our
case have profound consequences. They reduce the problem so that
on has effectively to diagonalize subspaces whose dimensions
increase polynomially with $r$ (rather than according to the power
law $3^r$). To formulate these features conveniently we introduce
the following conventions for state-labels.

For the fundamental case, $r=1$, the 3-dimensional basis is
denoted as
\begin{equation}
\left|1\right\rangle\equiv \left|\begin{array}{c} 1 \\ 0 \\  0
\\\end{array}\right\rangle,\,\,\,\left|2\right\rangle\equiv \left|\begin{array}{c}
   0 \\  1 \\  0 \\\end{array}\right\rangle,\,\,\,\left|3\right\rangle\equiv \left|\begin{array}{c}
0 \\  0 \\ 1 \\
\end{array}\right\rangle.
\end{equation}
For $r>1$, the order of the indices (1,2,3) represents the
tensored structure. Thus, for example, for $r=5$,
\begin{equation}
\left|11231\right\rangle \equiv
\left|1\right\rangle\otimes\left|1\right\rangle\otimes\left|2\right\rangle\otimes
\left|3\right\rangle\otimes\left|1\right\rangle.
\end{equation}
The fundamental realizations $t_{ij}^{(1)}$ of (2.2) implemented
in the tensored structure (1.14) of $t_{ij}^{(r)}$ lead to the
following major consequences.

\begin{description}
    \item[(I):] Each set of states corresponding to given sum of the indices (state labels 1,2,3) forms a closed
    subspace under the action of $T^{(r)}\left(\theta\right)$.
    Define with (with $a_i=\left(1,2,3\right)$)
    \begin{equation}
W\left|a_1a_2\ldots a_r\right\rangle=\left(a_1+a_2+\cdots+
a_r\right)\left|a_1a_2\ldots a_r\right\rangle.
\end{equation}
Then
\begin{equation}
\left[T^{(r)}\left(\theta\right),W\right]=0
\end{equation}
implying for each state on the right of
 \begin{equation}
T^{(r)}\left(\theta\right)\left|a_1a_2\ldots
a_r\right\rangle=\sum_{b_i}f_{(a,b)}\left(\theta\right)\left|b_1b_2\ldots
b_r\right\rangle,\qquad \left(b_1+b_2+\cdots+
b_r\right)=\left(a_1+a_2+\cdots+ a_r\right).
\end{equation}
Thus the $3^r$ dimensional base space of $T^{(r)}$ splits into
$(2r+1)$ closed subspaces under the action of $T^{(r)}$ as
\begin{equation}
S_r,\,S_{r+1},\ldots,\,S_{2r-1},\,S_{2r},\,S_{2r+1},\ldots,\,S_{3r},
\end{equation}
where $S_n$ corresponds to $a_1+a_2+\cdots+a_r=n$. In constructing
eigenstates of $T^{(r)}$ each $S_n$ can be treated separately
simplifying the problem considerably. The simplest subspaces are
the extreme ones, namely
\begin{equation}
S_r=\left|11\ldots 1\right\rangle
\end{equation}
and
\begin{equation}
S_{3r}=\left|33\ldots 3\right\rangle
\end{equation}
(the index 1(3) being repeated $r$ times). These are already
automatically eigenstates. The highest dimensional subspace is
obtained for $n=2r$ which includes the state $\left|22\ldots
2\right\rangle$. Special feature of some subspaces will be
displayed below.
    \item[(II):] Within each subspace again
    $T^{(r)}\left(\theta\right)$ commutes with circular
    permutations of states labels. Thus (CP) representing a
    circular permutation,
\begin{equation}
\left[T^{(r)}\left(\theta\right),(CP)\right]=0
\end{equation}
in the sense
\begin{equation}
(CP)^2T^{(r)}\left|a_1a_2a_3\ldots
a_{r-1}a_r\right\rangle=(CP)T^{(r)}\left|a_ra_1a_2\ldots
a_{r-2}a_{r-1}\right\rangle=T^{(r)}\left|a_{r-1}a_ra_1\ldots
a_{r-3} a_{r-2}\right\rangle
\end{equation}
and so on for all successive (CP) of the indices
$\left(a_1a_2\ldots a_r\right)$.

    \item[(III):] As a consequence the states in each invariant
    subspace can again be grouped together implementing roots of unity as follows. Let $\omega$ be any $r$-th root
    of unity, i.e.
\begin{equation}
\omega=\left(1,e^{i\frac{2\pi}r},e^{i\frac{2\pi}r\cdot
2},\ldots,e^{i\frac{2\pi}r\cdot (r-1)}\right)
\end{equation}
and (for each possible value of $\omega$, separately)
\begin{eqnarray}
&&\left|a_1a_2a_3\ldots a_{r-1}a_r\right\rangle_\omega\equiv
\left|a_1a_2a_3\ldots a_{r-1}a_r\right\rangle +\omega
\left|a_ra_1a_2\ldots a_{r-2}a_{r-1}\right\rangle+\omega^2
\left|a_{r-1}a_ra_1\ldots
a_{r-3}a_{r-2}\right\rangle\nonumber\\
&&\phantom{\left|a_1a_2a_3\ldots
a_{r-1}a_r\right\rangle_\omega\equiv}+\cdots+\omega^{r-1}
\left|a_2a_3a_4\ldots a_{r}a_{1}\right\rangle
\end{eqnarray}
The components states are, evidently, all in the same invariant
subspace. For $r$ different values of $\omega$ these provide a
mutually orthogonal set of $r$ states diagonalizing (CP) since
\begin{equation}
(CP)\left|a_1a_2a_3\ldots a_{r-1}a_r\right\rangle_\omega=\omega
\left|a_1a_2a_3\ldots a_{r-1}a_r\right\rangle_\omega
\end{equation}
The action of $T^{(r)}$ on, say, $\left|a_1a_2\ldots
a_{r-1}a_r\right\rangle$ gives directly, due to (3.10), that on
$\left|a_1a_2\ldots a_{r-1}a_r\right\rangle_\omega$ for all values
of $\omega$. Thus one can effectively reduce the dimension of the
relevant subspace $S_n$ for a given sum of state labels,
$\left(a_1+a_2+\cdots +a_{r-1}+a_r\right)=n$. Such a "two-step
reduction", firstly restriction to invariant subspaces $S_n$,
secondly introduction of roots of unity to form eigenstates of
(CP) will be shown to lead to a much slower increase with $r$ (as
compared to $e^{(\ln 3)r}$) of the dimension of the spaces on
which one has to diagonalize $T^{(r)}$. This will be first
displayed through particular examples. The general formulation
will be given at the end of this section.

\item[(IV):] But another symmetry is appropriately mentioned at
this stage (to be illustrated later explicitly). Interchanging the
indices as
\begin{equation}
\left(1,2,3\right)\longrightarrow \left(3,2,1\right)
\end{equation}
The action of $T^{(r)}$ is directly obtained via the inversion
\begin{equation}
q\longrightarrow q^{-1}
\end{equation}
in each coefficient. Thus the invariant subspaces related through
(3.14) need not ne studied separately. The corresponding
eigenstates and eigenvalues are related through (3.15). It is
sufficient to study the first $\left(r+1\right)$ subspaces since
under (3.14) and (3.15),
\begin{equation}
S_{2r}\longrightarrow S_{2r},\qquad
\left(S_r,S_{r+1},\ldots,S_{2r-1}\right) \rightleftharpoons
\left(S_{3r},S_{3r-1},\ldots,S_{2r+1}\right).
\end{equation}
Explicit examples for $r=\left(3,4\right)$ will follow. Our
$\hat{o}_N$ braid matrices remain nontrivial for $q=1$ as pointed
out in Ref. 2. Now (3.14) becomes a full symmetry. The degeneracy
thus induced is of interest.

\item[(V):] A final crucial feature is due to (1.16),
\begin{equation}
\left[T^{(r)}\left(\theta\right),T^{(r)}\left(\theta'\right)\right]=0.
\end{equation}
Suppose that for, say, $r=4$ in some subspace one obtains a closed
subset of states $\left(A,B,C,D\right)$ with
\begin{equation}
T^{(4)}\left(\theta\right)A=a_{11}A+a_{12}B+a_{13}C+a_{14}D,\ldots,\,
T^{(4)}\left(\theta\right)D=d_{11}A+d_{12}B+d_{13}C+d_{14}D.
\end{equation}
The coefficients $\left(a_{11},\ldots,d_{14}\right)$ are in
general polynomials in $K\left(\theta\right)$, the maximal degree
being $r=4$ for this case. Define eigenstates as
\begin{equation}
T^{(4)}\left(\theta\right)\left(\alpha A+\beta B+\gamma C+\delta
D\right)=v \left(\alpha A+\beta B+\gamma C+\delta D\right)
\end{equation}
which are to be solved for by implementing (3.18) on the left.
Consistency with (3.17) imposes $\theta$-independence
($K$-independence) of $\left(\alpha,\beta,\gamma,\delta\right)$.
Hence on the right only $v$ can be $K$-dependent. All
$K$-dependence of $\left(a_{11},\ldots,d_{14}\right)$ on the left
must thus factorize as a polynomial (here for $r=4$)
\begin{equation}
v=f_4 K^4+f_3 K^3+f_2 K^2+f_1 K+f_0
\end{equation}
for suitable $\left(f_4,f_3,f_2,f_1,f_0\right)$ which can depend
on $\left(q,\omega\right)$ only. In general this leads to a set of
overdetermined set of coupled linear equations (for our case) in
\begin{equation}
\left(\alpha,\beta,\gamma,\delta; f_4,f_3,f_2,f_1,f_0\right)
\end{equation}
Varied illustrations will follows. Moreover, while all eigenvalues
are, in general, $K$- and $q$- dependent, all explicit
$q$-dependence (except for the implicit one through $K$ of (1.17)
and (1.18)) must cancel in the overall trace (summing over all
subspaces) to give (2.9), i.e.
\begin{equation}
\hbox{Tr}\left(T^{(r)}\right)=3\left(1+K\right)^r.
\end{equation}
This provides stringent check (Appendix A).
\end{description}

\paragraph{\large\bf Special features of the subspaces
$\left(S_{r},S_{3r}\right)$, $\left(S_{r+1},S_{3r-1}\right)$,
$S_{2r}$:}

\paragraph{\bf $\bullet$ \underline{$\left(S_{r},S_{3r}\right)$}:} As
mentioned following (3.7), (3.8) these two are 1-dimensional
subspaces. One obtains immediately, for all $r$,
\begin{eqnarray}
&&T^{(r)}\left(\theta\right)\left|11\ldots 1\right\rangle=\left(1+K^r\right)\left|11\ldots 1\right\rangle\\
&&T^{(r)}\left(\theta\right)\left|33\ldots
3\right\rangle=\left(1+K^r\right)\left|33\ldots 3\right\rangle
\end{eqnarray}
These eigenstates of (CP), singlets, provide the simplest
illustrations of (3.14), (3.15).

\paragraph{\bf $\bullet$ \underline{$\left(S_{r+1},S_{3r-1}\right)$}:} For arbitrary $r$,
with
\begin{equation}
\omega= \left(1,e^{i\frac{2\pi}r},e^{i\frac{2\pi}r\cdot
2},\ldots,e^{i\frac{2\pi}r\cdot (r-1)}\right)
\end{equation}
define
\begin{eqnarray}
&&X_\omega=\left|111\ldots 12\right\rangle+\omega\left|211\ldots
11\right\rangle+\omega^2\left|121\ldots
11\right\rangle+\omega^{r-1}\left|111\ldots
21\right\rangle,\\
&&Y_\omega=\left|333\ldots 32\right\rangle+\omega\left|233\ldots
33\right\rangle+\omega^2\left|323\ldots
33\right\rangle+\omega^{r-1}\left|333\ldots 23\right\rangle.
\end{eqnarray}
One easily obtains
\begin{eqnarray}
&&T^{(r)}\left(\theta\right)X_\omega=\left(K^r\omega+\omega^{r-1}\right)X_\omega,\\
&&T^{(r)}\left(\theta\right)Y_\omega=\left(K^r\omega+\omega^{r-1}\right)Y_\omega.
\end{eqnarray}
For the $r$ values of $\omega$ one obtains thus, in a single
stroke, all the requisite $r$ eigenstates for these two
$r$-dimensional subspaces. Note that
\begin{equation}
\sum_\omega\left(K^r\omega+\omega^{r-1}\right)=\sum_\omega\left(K^r\omega+\omega^{-1}\right)=0.
\end{equation}
Hence $\left(S_{r+1},S_{3r-1}\right)$ do not contribute to the
total trace $Tr\left(T^{(r)}\left(\theta\right)\right)$.

For $S_{r+2}$ $\left(S_{3r-2}\right)$ already the structure of
eigenstates and eigenvalues are not so simple. (See App. A for
$r=3,4$). Some special features of $S_{2r}$ are however worth
mentioning, particularly to compare the structures of $r$ prime
and non-prime.

\paragraph{\bf $\bullet$ \underline{$\left(S_{2r}\right)$}:} Like $\left|11\ldots
1\right\rangle$ and $\left|33\ldots 3\right\rangle$,
$\left|22\ldots 2\right\rangle$ is also a singlet under (CP). But
unlike the former the latter one does not form an 1-dimensional
subspace. It can get coupled with the other states of $S_{2r}$
(for $\omega=1$) as follows. When $r$ is prime, apart form
$\left|22\ldots 2\right\rangle$, $S_{2r}$ is composed of $r$-plets
(formed using $\omega$ with $\omega^r=1$). When  $r$ is
factorizable there can be intermediate multiplets corresponding to
factors $\left(n_1,n_2,\ldots,n_k\right)$ of $r=\left(n_1n_2\cdots
n_k\right)$. Thus for $r=4$ (the first factorizable $r$) there are
doublets corresponding to $r=2\cdot 2$. For $r=6$, there are
doublets and triplets between $1-$ and $6-$plets. Let us
illustrate the situation using the simplest non-trivial cases
$r=3,4$.

\paragraph{\bf $\diamond$ \underline{$\left(r=3,S_{6}\right)$}:} Define
\begin{equation}
A_1=\left|222\right\rangle,\qquad
B_\omega=\left|123\right\rangle+\omega\left|312\right\rangle+\omega^2\left|231\right\rangle,\qquad
C_\omega=\left|321\right\rangle+\omega\left|132\right\rangle+\omega^2\left|213\right\rangle,
\end{equation}
where $\omega=\left(1,e^{i\frac{2\pi}3},e^{i\frac{2\pi}3\cdot
2}\right)$. In our notation $A_1$ indicates that here (for
singlet) one has only $\omega=1$. Correspondingly
$\left(B_1,C_1\right)$ will denote the latter for $\omega=1$.
Consistently with (3.9) set
\begin{equation}
T^{(3)}\left(\theta\right)\left(\alpha A_1+\beta B_1+\gamma
C_1\right)=v\left(\alpha A_1+\beta B_1+\gamma C_1\right)
\end{equation}
for $\omega=1$ and
\begin{equation}
T^{(3)}\left(\theta\right)\left(\mu B_\omega+\nu
C_\omega\right)=w\left(\mu B_\omega+\nu C_\omega\right)
\end{equation}
for $\omega=\left(e^{i\frac{2\pi}3},e^{-i\frac{2\pi}3}\right)$.
Here $\left(v,w\right)$ are assumed to be cubic polynomials in $K$
and $\left(\alpha,\beta,\gamma\right)$, $\left(\mu,\nu\right)$ to
be $K$-independent. Note also that
\begin{equation}
\left(1\rightleftharpoons 3\right) B_\omega=C_\omega.
\end{equation}
Hence (consistently with (3.14), (3.15)) one obtains the
coefficients in $T^{(3)}\left(\theta\right)C_\omega$ by inverting
$q$ to $q^{-1}$ in those of $T^{(3)}\left(\theta\right)B_\omega$.
Explicit solutions are given in App. A. Here we only note that the
decoupling of $A_1$ in (3.33) is assured via the structure
\begin{eqnarray}
&&T^{(3)}\left(\theta\right)A_1=a_{11}A_1+a_{12}B_1+a_{13}C_1\nonumber\\
&&T^{(3)}\left(\theta\right)B_\omega=\left(1+\omega+\omega^2\right)b_{11}A_1+b_{12}B_\omega+b_{13}C_\omega\nonumber\\
&&T^{(3)}\left(\theta\right)C_\omega=\left(1+\omega+\omega^2\right)c_{11}A_1+c_{12}B_\omega+c_{13}C_\omega,
\end{eqnarray}

\paragraph{\bf $\diamond$ \underline{$\left(r=4,S_{8}\right)$}:} Here, after the
$\left(CP\right)$-singlet
\begin{equation}
A_1=\left|2222\right\rangle
\end{equation}
one has also the doublets
\begin{equation}
B_{\pm 1}=\left|1313\right\rangle\pm\left|3131\right\rangle
\end{equation}
and then the quartets completing the 19 dimensional $S_8$ for all
values of $\omega$, namely,
\begin{eqnarray}
&&\omega= \left(1,e^{i\frac{2\pi}4},e^{i\frac{2\pi}4\cdot
2},e^{i\frac{2\pi}4\cdot 3}\right)=\left(1,i,-1,-i\right)\\
&&C_\omega=\left|1133\right\rangle+\omega\left|3113\right\rangle+\omega^2\left|3311\right\rangle+
\omega^3\left|1331\right\rangle,\nonumber\\
&&D_\omega=\left|1223\right\rangle+\omega\left|3122\right\rangle+\omega^2\left|2312\right\rangle+
\omega^3\left|2231\right\rangle,\nonumber\\
&&E_\omega=\left|3221\right\rangle+\omega\left|1322\right\rangle+\omega^2\left|2132\right\rangle+
\omega^3\left|2213\right\rangle,\nonumber\\
&&F_\omega=\left|1232\right\rangle+\omega\left|2123\right\rangle+\omega^2\left|3212\right\rangle+
\omega^3\left|2321\right\rangle.
\end{eqnarray}
Note also that
\begin{equation}
\left(1\rightleftharpoons 3\right)
\left(C_\omega,D_\omega,E_\omega,F_\omega\right)=\left(\omega^2
C_\omega,E_\omega,D_\omega,\omega^2 F_\omega\right)
\end{equation}
which simplifies computations according to (3.14), (3.15). The set
$F_\omega$ alone has a distinctive feature. The two indices 2
remain separated (unlike for $D_\omega$, $E_\omega$) under (CP).
This singles it out as directly an eigenstates of
$T^{(4)}\left(\theta\right)$ (App. A). As for $\left(C_\omega,
D_\omega,E_\omega\right)$ decouplings, analogous to (3.35) but in
two stages
\begin{description}
    \item[(1)] from $A_1$ for $\omega=\left(-1,\pm i\right)$
    \item[(2)] and also from $B_{\pm 1}$ for $\omega=\left(\pm i\right)$
\end{description}
are assured through factors of the type (App. A)
\begin{equation}
\left(1+\omega+\omega^2+\omega^3\right),\,
\left(1+\omega^2\right)\left(1\pm \omega\right).
\end{equation}
The maximal set of 5 coupled linear equations arises for
($\omega=1$)
\begin{equation}
T^{(4)}\left(\theta\right)\left(aA_1+bB_1+cC_1+dD_1+eE_1\right)=v_1\left(aA_1+bB_1+cC_1+dD_1+eE_1\right).
\end{equation}
To conclude we emphasize again that for $r=(3,4)$ in base spaces
respectively of dimensions (27,81) the maximal set of coupled
linear equations encountered are sets of (3,5) respectively. This
is the slow growth with $r$ signalled before (end of (III)).

For $r=\left(1,2,3,4\right)$ we have studied the invariant
subspaces $S_n$ explicitly. Let us now indicate the general
situation. Associate the variables $\left(x_1,x_2,x_3\right)$ to
the states
$\left(\left|1\right\rangle,\left|2\right\rangle,\left|3\right\rangle\right)$
respectively. In the expansion
\begin{equation}
\left(x_1+x_2+x_3\right)^r=\sum_{n_1,n_2,n_3}C_{n_1,n_2,n_3}x_1^{n_1}x_2^{n_2}x_3^{n_3}.
\end{equation}
for each term
\begin{equation}
n_1+n_2+n_3=r
\end{equation}
and
\begin{equation}
\sum_{n_1,n_2,n_3}C_{n_1,n_2,n_3}=3^r.
\end{equation}
Imposing an additional constraint one obtains the subsets
\begin{equation}
\dim S_n=\sum_{n_1,n_2,n_3}C_{n_1,n_2,n_3},\qquad
\left(n_1+2n_2+3n_3=n\right)
\end{equation}
for $n=\left(r,r+1,\ldots,2r,\ldots,3r\right)$. The dimension of
the total base space for order $r$ is given by (3.45).

Let us consider as an example the central subspace $S_{10}$ for
$r=5$. From (3.46) one easily finds
\begin{equation}
\dim S_{10}=51\qquad \left(r=5\right).
\end{equation}
The states can be grouped into multiplets as follows with
\begin{eqnarray}
&&\omega^5=1,\\
&&V_1^{(1)}=\left|22222\right\rangle,\nonumber\\
&&V_2^{(\omega)}=\left(\left|13222\right\rangle+\omega \left|21322\right\rangle+\omega^2
\left|22132\right\rangle+\omega^3\left|22213\right\rangle+\omega^4\left|32221\right\rangle\right)\equiv
\left((CP)\left|13222\right\rangle\right)_\omega,\nonumber\\
&&V_3^{(\omega)}=\left((CP)\left|12322\right\rangle\right)_\omega,\nonumber\\
&&V_4^{(\omega)}=\left((CP)\left|13132\right\rangle\right)_\omega,\nonumber\\
&&V_5^{(\omega)}=\left((CP)\left|13312\right\rangle\right)_\omega,\nonumber\\
&&V_6^{(\omega)}=\left((CP)\left|13321\right\rangle\right)_\omega,\nonumber\\
&&\left(V_7^{(\omega)},\ldots,V_{11}^{(\omega)}\right)=\left(1\rightleftharpoons
3\right) \left(V_2^{(\omega)},\ldots,V_{6}^{(\omega)}\right),
\end{eqnarray}
i.e.
$V_7^{(\omega)}=\left((CP)\left|31222\right\rangle\right)_\omega$
and so on. For $\omega=1$, now one has to solve a set of 11
(coupling $V_1^{(1)},\ldots V_{11}^{(1)}$) linear equations. This
is the maximal such set for $r=5$ where the total dimension is
243.

Whenever $r$ is a prime number, i.e.
$r=\left(1,3,5,7,11,13,\ldots\right)$, the multiplet structure is
relatively simple. Thus for $S_{2r}$ apart from $\left|22\ldots
2\right\rangle$ there are only $r$-plets in terms of the roots
$\omega^r=1$. When $r$ is factorizable lower multiplets can arise
corresponding to factors of $r$. We have illustrated this for
$r=4$.

\section{Chain Hamiltonians ($N=3$)}
\setcounter{equation}{0}

The Hamiltonian for order $r$ is defined as
\begin{equation}
H^{(r)}=\left(T^{(r)}\left(\theta\right)\right)^{-1}_{\theta=0}\left(\partial_\theta
\left(T^{(r)}\left(\theta\right)\right)\right)_{\theta=0}
\end{equation}
Instead of using the standard formulation as a sum (see the basic
references in sec. 4 of Ref. 3)
\begin{equation}
H^{(r)}=\sum_{k=1}^rI\otimes I\otimes\cdots\otimes
\dot{\hat{R}}_{k,k+1}\left(0\right)\otimes I\otimes\cdots\otimes
I,
\end{equation}
where
\begin{equation}
\dot{\widehat{R}}_{k,k+1}\left(0\right)=\left(\partial_\theta
\hat{R}_{k,k+1}\left(\theta\right)\right)_{\theta=0}
\end{equation}
with the circular boundary condition for $k=r$ ($r+1\approx 1$) we
we will use (4.1) directly, as explained below, in a fashion
particularly well-adapted to our formalism for constructing
eigenstates.

Define starting from (1.17) i.e.
\begin{eqnarray}
&&K\left(\theta\right)=-\frac{\sinh
\theta}{\sinh\left(\eta+\theta\right)},\\
&&\dot{K}_0\equiv \left(\partial_\theta
K\left(\theta\right)\right)_{\theta=0}=-\left(\sinh\eta\right)^{-1}
\end{eqnarray}
with $K_0=\left(K\left(\theta\right)\right)_{\theta=0}=0$. We
start with eigenstate of $T^{(r)}\left(\theta\right)$
\begin{equation}
\left|V\right\rangle_\omega=\left(c_1A_1+c_2A_2+\cdots+c_mA_m\right)_\omega,
\end{equation}
where the subscript $\omega$ indicates that each $A_i$
($i=1,\ldots,m$) is an eigenstate of (CP), circular permutation of
$r$ state labels corresponding to a subspace $S_n$
($n=r,\ldots,3r$). (See Sec. 3 and App. A). Thus for example, for
$r=3$ and $S_n=S_5$ (see  (A.13) and (A.18)-(A.20))
$\left|V\right\rangle_\omega=aA_\omega+bB_\omega$, where
\begin{equation}
A_\omega=\left(\left|113\right\rangle+\omega\left|311\right\rangle+\omega^2\left|131\right\rangle\right),\qquad
B_\omega=\left(\left|122\right\rangle+\omega\left|212\right\rangle+\omega^2\left|221\right\rangle\right),
\end{equation}
with $\omega^3=1$. Quite generally, if for (4.6)
\begin{equation}
T^{(r)}\left(\theta\right)\left|V\right\rangle=v\left|V\right\rangle=v\left(\sum_kc_kA_k\right)
\end{equation}
then as explained and emphasized (in sec. 3 and App. A) the
coefficients $c_k$ can depend on $q$ (but not on $\theta$) the
only $\theta$-dependence on the right is in $v$, a polynomial of
order $r$ in $K\left(\theta\right)$,
\begin{equation}
v=f_r\left(K\left(\theta\right)\right)^r+f_{r-1}\left(K\left(\theta\right)\right)^{r-1}+\cdots+
f_1\left(K\left(\theta\right)\right)+f_0,
\end{equation}
where the coefficients $f_i$ are each $\theta$-independent. Thus
for (4.7) the solutions (for each value of $\omega$) are
\begin{eqnarray} &&\hbox{\bf (1)}\qquad
\left(a,b\right)=\left(q^{1/2}+\omega
q^{-1/2},1\right),\\
&&\phantom{\hbox{\bf (1)}\qquad} v=\omega^2K^3+
\left(\left(q+q^{-1}\right)\omega^2+\left(1+\omega+\omega^2\right)\right)K^2+
\left(\left(q+q^{-1}\right)\omega+\left(1+\omega+\omega^2\right)\right)K+\omega\nonumber\\
&&\hbox{\bf (2)}\qquad
\left(a,b\right)=\left(1,-\left(q^{1/2}+\omega^2
q^{-1/2}\right)\right),\qquad v=\omega^2K^3+\omega.
\end{eqnarray}
From (4.5), (4.8) and (4.9) one obtains (since $K_0=0$) the
general result (with $\dot{T}_0^{(r)}\equiv \left(\partial_\theta
T^{(r)}\left(\theta\right)\right)_{\theta=0}$, $T_0^{(r)}=\left(
T^{(r)}\left(\theta\right)\right)_{\theta=0}$)
\begin{eqnarray}
&&\dot{T}_0^{(r)}\left|V\right\rangle=\dot{K}_0f_1\left|V\right\rangle,\\
&&T_0^{(r)}\left|V\right\rangle=f_0\left|V\right\rangle=\omega\left|V\right\rangle
\end{eqnarray}
and hence
\begin{equation}
\left(T_0^{(r)}\right)^{-1}\left|V\right\rangle=\omega^{r-1}\left|V\right\rangle
\end{equation}
The result $f_0=\omega$ (and $f_0^{-1}=\omega^{r-1}$ for
$\omega^r=1$) is a general one. This corresponds to our use
eigenstates of (CP) as basis states since for our class $T_0$
coincides with (CP). Hence finally
\begin{equation}
H^{(r)}\left|V\right\rangle=T_0^{-1}\dot{T}_0\left|V\right\rangle=
\left(\dot{K}_0\omega^{r-1}f_1\right)\left|V\right\rangle
\end{equation}
Thus starting with an eigenstate of $T^{(r)}\left(\theta\right)$
in our formalism it remains one of $H^{(r)}$ and the eigenvalue of
$H^{(r)}$ is extracted, as above from that of $T^{(r)}$. Note that
for
\begin{equation}
f_1=0, \qquad H^{(r)}\left|V\right\rangle=0.
\end{equation}
Thus for (4.11)
\begin{equation}
H^{(3)}\left(A_\omega-\left(q^{1/2}+\omega^2q^{-1/2}\right)B_\omega\right)=0.
\end{equation}
From (3.23)-(3.29) it follows that, for all $r$,
\begin{equation}
H^{(r)}\left(S_r,S_{3r};S_{r+1},S_{3r-1}\right)\approx 0,
\end{equation}
i.e. each eigenstate belonging to these subspaces is annihilated
by $H^{(r)}$.

For $r=2$, the explicit form of the Hamiltonian is
\begin{eqnarray}
&&\left(\dot{K}_0\right)^{-1}H^{(2)}=\left(q+q^{-1}\right)\left(11\right)\otimes (33)+
\left(q^{1/2}+q^{-1/2}\right)\left(12\right)\otimes (32)+2\left(13\right)\otimes (31)\nonumber\\
&&\phantom{\left(\dot{K}_0\right)^{-1}H^{(2)}=}+\left(q^{1/2}+q^{-1/2}\right)\left(21\right)\otimes
(23)+2\left(22\right)\otimes
(22)+\left(q^{1/2}+q^{-1/2}\right)\left(23\right)\otimes (21)\nonumber\\
&&\phantom{\left(\dot{K}_0\right)^{-1}H^{(2)}=}+2\left(31\right)\otimes
(33)+\left(q^{1/2}+q^{-1/2}\right)\left(32\right)\otimes
(12)+\left(q+q^{-1}\right)\left(33\right)\otimes (11)
\end{eqnarray}
Consistently with (4.18)
\begin{equation}
H^{(2)}\left(\left|11\right\rangle,\left|33\right\rangle;\left|12\right\rangle,\left|21\right\rangle;
\left|23\right\rangle,\left|32\right\rangle\right)=0.
\end{equation}
For the only remaining subspace $S_4$, setting
\begin{equation}
H^{(2)}\left(a\left|13\right\rangle+b\left|22\right\rangle+c\left|31\right\rangle\right)=
v_H\left(a\left|13\right\rangle+b
\left|22\right\rangle+c\left|31\right\rangle\right).
\end{equation}
One obtains the solutions
\begin{eqnarray}
&&\hbox{\bf (1)}\qquad
\left(a,b,c\right)=\left(1,-\left(q^{1/2}+q^{-1/2}
\right),1\right);\,\,v_H=0,\nonumber\\
&&\hbox{\bf (2)}\qquad \left(a,b,c\right)=\left(1,0,-1\right);\,\,v_H=\dot{K}_0\left(q+q^{-1}-2\right),\nonumber\\
&&\hbox{\bf (3)}\qquad
\left(a,b,c\right)=\left(1,2\left(q^{1/2}+q^{-1/2}\right)^{-1},1\right);\,\,
v_H=\dot{K}_0\left(q+q^{-1}+4\right).
\end{eqnarray}
Combining (4.15) with (A.4)-(A.6) one consistently reproduces the
results (4.20)-(4.22) obtained using the explicit form (4.19).
Note that for $r=2$ and $\omega=-1$, (4.15) gives
\begin{equation}
H^{(2)}\left|V\right\rangle=-\dot{K}_0f_1\left|V\right\rangle.
\end{equation}
This corresponds to the positive sign in solution (2) of (4.22)
since in (A.6) the corresponding factor is
\begin{equation}
f_1=-\left(q+q^{-1}-2\right)
\end{equation}
Such changes of sign introduce a qualitative change:
$\hbox{Tr}\left(T^{(r)}\right)$ in (2.9) has no explicit
dependence on $q$ (only an implicit one though $K$). But
$\hbox{Tr}\left(H^{(r)}\right)$ can have explicit $q$-dependence.
For the simple example above $\left(N=3,r=2\right)$
\begin{equation}
\hbox{Tr}\left(H^{(2)}\right)=2\dot{K}_0\left(q+q^{-1}+1\right)
\end{equation}

\smallskip

\noindent{\large\bf Selection rules for transitions:} Adopting the
convention of attaching to the states $\left(\left|1\right\rangle,
\left|2\right\rangle,\right.$ $\left.\left|3\right\rangle\right)$
respectively the "spins"
\begin{equation}
\left(+,0,-\right)
\end{equation}
it is seen from (4.19) that the action of the Hamiltonian on
neighboring sites, induces transitions only when the sum of the
two spins is zero, i.e. for
\begin{equation}
\left(+-\right),\qquad \left(00\right),\qquad\left(-+\right)
\end{equation}
The final states corresponding again to zero sum. Thus one has
non-zero matrix elements for a neighboring pair
$\left|ij\right\rangle{H\atop\rightarrow} \left|kl\right\rangle$
only when for the corresponding spins
\begin{equation}
\sigma_i+\sigma_j=0=\sigma_k+\sigma_l
\end{equation}
Such matrix elements depend on $\left(\dot{K}_0,q\right)$. The
structure of $H^{(r)}$ in (4.2) indicates that (4.21) is a generic
feature. Any pair of the type (4.27) somewhere in the chain can
start transitions which can propagate along the chain since the
three possibilities in (4.27) can create such a pair with the next
neighboring site and so on.

\section{$N>3$}
\setcounter{equation}{0}

Three basic features displayed and studied at length for $N=3$
are:
\begin{description}
    \item[(1)] A simple recursion relation yielding the trace of
    the transfer matrix for any order $r$. (See (2.2)-(2.9)).
    \item[(2)] Invariant subspaces corresponding to the sum of the
    state labels. (See (3.1)-(3.6) and App. A.)
    \item[(3)] Role of (CP) circular permutation of state labels
    within each invariant subspace $S_n$. (See (3.9)-(3.13) and App.
    A.)
\end{description}
It was shown (for $N=3$) how (2) and (3) greatly simplify the
construction of eigenstates and eigenvalues of
$T^{(r)}\left(\theta\right)$ for successive values of $r$.

We now indicate how these features are carried over for $N>3$ via
the simplest possibilities, namely $N=4$, $r=\left(1,2\right)$.
Now, as compared to (1.17)-(1.18)
\begin{equation}
K\left(\theta\right)=-\frac{\sinh\theta}{\sinh\left(\eta+\theta\right)}
\end{equation}
where $e^\eta+e^{-\eta}=\left(q^2+1+q^{-2}\right)+1
=\left(q+q^{-1}\right)^2$. As compared to (2.2) (writing $t_{ij}$
for $t_{ij}^{(1)}\left(\theta\right)$, $i=\left(1,2,3,4\right)$
and $K$ for $K\left(\theta\right)$)
$t_{11}=\left|\begin{array}{cccc}
   1 & 0 & 0 & 0 \\
   0 & 0 & 0 & 0\\
   0 & 0 & 0 & 0\\
   0 & 0 & 0 & K\\
\end{array}\right|\equiv \left(11\right)+K\left(44\right)$ and
similarly,
\begin{eqnarray}
&&t_{12}=\left(21\right)+Kq\left(43\right),\,\,\,
t_{13}=\left(31\right)+Kq\left(42\right),\,\,\,
t_{14}=\left(1+Kq^2\right)\left(41\right),\nonumber\\
&&t_{21}=\left(12\right)+Kq^{-1}\left(34\right),\,\,\,
t_{22}=\left(22\right)+K\left(33\right),\,\,\,
t_{23}=\left(1+K\right)\left(32\right),\,\,\,
t_{24}=Kq\left(31\right)+\left(42\right),\nonumber\\
&&t_{31}=\left(13\right)+Kq^{-1}\left(24\right),\,\,\,
t_{32}=\left(1+K\right)\left(23\right),\,\,\,
t_{33}=K\left(22\right)+\left(33\right),\,\,\,
t_{34}=Kq\left(21\right)+\left(43\right),\\
&&t_{41}=\left(1+Kq^{-2}\right)\left(14\right),\,\,\,
t_{42}=Kq^{-1}\left(13\right)+\left(24\right),\,\,\,
t_{43}=Kq^{-1}\left(12\right)+\left(34\right),\,\,\,
t_{44}=K\left(11\right)+\left(44\right)\nonumber
\end{eqnarray}
As compared to (2.3)-(2.9) recursion relations are now
(suppressing arguments $\theta$)
\begin{equation}
t_{ij}^{(r+1)}=t_{i1}^{(1)}\otimes
t_{1j}^{(r)}+t_{i2}^{(2)}\otimes t_{2j}^{(r)}+t_{i3}^{(1)}\otimes
t_{3j}^{(r)}+t_{i4}^{(1)}\otimes t_{4j}^{(r)},\qquad
\left(i,j=1,2,3,4\right)
\end{equation}
giving (due to (5.2))
\begin{eqnarray}
&&t_{1j}^{(r+1)}=\left(\left(11\right)+K\left(44\right)\right)\otimes
t_{1j}^{(r)}+\left(\left(21\right)+Kq\left(43\right)\right)\otimes
t_{2j}^{(r)}+\left(\left(31\right)+Kq\left(42\right)\right)\otimes
t_{3j}^{(r)}\nonumber\\
&&\phantom{t_{1j}^{(r+1)}=}
+\left(1+Kq^2\right)\left(41\right)\otimes t_{4j}^{(r)},\nonumber\\
&&t_{2j}^{(r+1)}=\left(\left(12\right)+Kq^{-1}\left(34\right)\right)\otimes
t_{1j}^{(r)}+\left(\left(22\right)+K\left(33\right)\right)\otimes
t_{2j}^{(r)}+\left(1+K\right)\left(32\right)\otimes
t_{3j}^{(r)}
\nonumber\\
&&\phantom{t_{1j}^{(r+1)}=}
+\left(Kq\left(31\right)+\left(42\right)\right)\otimes t_{4j}^{(r)},\nonumber\\
&&t_{3j}^{(r+1)}=\left(\left(13\right)+Kq^{-1}\left(24\right)\right)\otimes
t_{1j}^{(r)}+\left(1+K\right)\left(23\right)\otimes
t_{2j}^{(r)}+\left(K\left(22\right)+\left(33\right)\right)\otimes
t_{3j}^{(r)}\nonumber\\
&&\phantom{t_{1j}^{(r+1)}=}+\left(Kq\left(21\right)+K\left(43\right)\right)\otimes t_{4j}^{(r)},\nonumber\\
&&t_{4j}^{(r+1)}=\left(1+Kq^{-2}\right)\left(14\right)\otimes
t_{1j}^{(r)}+\left(Kq^{-1}\left(13\right)+\left(24\right)\right)\otimes
t_{2j}^{(r)}+\left(Kq^{-1}\left(12\right)+\left(34\right)\right)\otimes
t_{3j}^{(r)}\nonumber\\
&&\phantom{t_{1j}^{(r+1)}=}+\left(K\left(11\right)+\left(44\right)\right)\otimes
t_{4j}^{(r)}.
\end{eqnarray}
Hence for the transfer matrix
\begin{eqnarray}
&&T^{(r+1)}=t_{11}^{(r+1)}+t_{22}^{(r+1)}+t_{33}^{(r+1)}+t_{44}^{(r+1)}\nonumber\\
&&\phantom{T^{(r+1)}}=\left|\begin{array}{cccc}
   t_{11}^{(r)}+Kt_{44}^{(r)} & t_{12}^{(r)}+Kq^{-1}t_{34}^{(r)} & t_{13}^{(r)}+Kq^{-1}t_{34}^{(r)}
   & \left(1+Kq^{-2}\right)t_{14}^{(r)} \\
    t_{21}^{(r)}+Kt_{43}^{(r)} & t_{22}^{(r)}+Kt_{33}^{(r)} & \left(1+K\right)t_{23}^{(r)}
   & Kq^{-1}t_{13}^{(r)}+t_{24}^{(r)} \\
    t_{31}^{(r)}+Kqt_{42}^{(r)} & \left(1+K\right)t_{32}^{(r)} & Kt_{22}^{(r)}+t_{33}^{(r)}
   & Kq^{-1}t_{12}^{(r)}+t_{34}^{(r)} \\
    \left(1+Kq^2\right)t_{41}^{(r)} & t_{42}^{(r)}+Kqt_{31}^{(r)} & Kqt_{21}^{(r)}+t_{43}^{(r)}
   & Kt_{11}^{(r)}+ t_{44}^{(r)}\\
\end{array}\right|.
\end{eqnarray}
One now obtains
\begin{equation}
\hbox{Tr}\left(T^{(r+1)}\right)=\left(K+1\right)\hbox{Tr}\left(T^{r}\right).
\end{equation}
But
\begin{equation}
T^{(1)}=\left(K+1\right)I_4.
\end{equation}
Hence
\begin{equation}
\hbox{Tr}\left(T^{(r)}\right)=4\left(K+1\right)^r.
\end{equation}
It is not difficult to obtain the general result (following from
the fact that only the diagonal blocks $t_{ii}^{(r)}$ have
diagonal terms)
\begin{equation}
\hbox{Tr}\left(T^{(r)}\right)=N\left(K+1\right)^r
\end{equation}
For $N=(3,4)$ the particular solutions are given by (2.9) and
(5.8) respectively.

Now let us consider the eigenstates of $T^{(r)}\left(\theta
\right)$ for $N=4$, $r=1,2$. As compared to (3.1) we now have the
fundamental state vectors
\begin{equation}
\left|1\right\rangle=\left|\begin{array}{c}1 \\0 \\ 0 \\0 \\
\end{array}\right\rangle,\,\left|2\right\rangle=\left|\begin{array}{c}0 \\1 \\ 0 \\0 \\
\end{array}\right\rangle,\left|3\right\rangle=\left|\begin{array}{c}0 \\0 \\ 1 \\0 \\
\end{array}\right\rangle,\left|4\right\rangle=\left|\begin{array}{c}0 \\0 \\ 0 \\1 \\
\end{array}\right\rangle
\end{equation}
and, as before, we denote tensor products as
\begin{equation}
\left|ijk\ldots\right\rangle=\left|i\right\rangle\otimes
\left|j\right\rangle\otimes \left|k\right\rangle\otimes\cdots
\end{equation}
For a given $n$, as before, the set of states with
\begin{equation}
i+j+k+\cdots=n
\end{equation}
will constitute the basis of the subspace $S_n$. For $r=1$ the
situation is trivial. From (5.2)
\begin{eqnarray}
&&T^{(1)}=t_{11}+t_{22}+t_{33}+t_{44}=\left(K+1\right)I_4,\\
&&T^{(1)}\left|i\right\rangle=\left(K+1\right)\left|i\right\rangle,\qquad
\left(i=1,2,3,4\right).
\end{eqnarray}
For $r=2$, from (5.2) and (5.5),
\begin{eqnarray}
&&T^{(2)}=\left(K^2+1\right)P+2K\left((11)\otimes (44)+(22)\otimes
(33)+(33)\otimes (22)+(44)\otimes (11)\right)\nonumber\\
&&\phantom{T^{(2)}=}K\left(q+q^{-1}\right)\left((12)\otimes
(43)+(13)\otimes (42)+(21)\otimes (34)+(24)\otimes
(31)+(31)\otimes (24)\right. \\
&&\phantom{T^{(2)}=}\left.+(34)\otimes (21)+(42)\otimes
(13)+(43)\otimes
(12)\right)+K\left(q^2+q^{-2}\right)\left((14)\otimes
(41)+(41)\otimes (14)\right),\nonumber
\end{eqnarray}
where
\begin{equation}
P=\sum_{ij}(ij)\otimes (ji),\qquad \left(i,j=1,2,3,4\right).
\end{equation}

Implementing the definitions (5.10)-(5.12) one obtains from (5.15)
for the subspaces $\left(S_2,\ldots,S_8\right)$ the following
results (with $\epsilon=\pm 1$)
\begin{eqnarray}
&&S_2:\qquad T^{(2)}\left|11\right\rangle=\left(K^2+1\right)
\left|11\right\rangle, \nonumber\\
&&S_8:\qquad T^{(2)}\left|44\right\rangle=\left(K^2+1\right)
\left|44\right\rangle, \nonumber\\
&&S_3:\qquad
T^{(2)}\left(\left|12\right\rangle+\epsilon\left|21\right\rangle\right)=\epsilon\left(K^2+1\right)
\left(\left|12\right\rangle+\epsilon\left|21\right\rangle\right), \nonumber\\
&&S_7:\qquad
T^{(2)}\left(\left|43\right\rangle+\epsilon\left|34\right\rangle\right)=\epsilon\left(K^2+1\right)
\left(\left|43\right\rangle+\epsilon\left|34\right\rangle\right), \nonumber\\
&&S_4:\qquad T^{(2)}\left|22\right\rangle=\left(K^2+1\right)
\left|22\right\rangle, \nonumber\\
&&\phantom{S_4:}\qquad
T^{(2)}\left(\left|13\right\rangle+\epsilon\left|31\right\rangle\right)=\epsilon\left(K^2+1\right)
\left(\left|13\right\rangle+\epsilon\left|31\right\rangle\right), \nonumber\\
&&S_6:\qquad T^{(2)}\left|33\right\rangle=\left(K^2+1\right)
\left|33\right\rangle, \nonumber\\
&&\phantom{S_4:}\qquad
T^{(2)}\left(\left|42\right\rangle+\epsilon\left|24\right\rangle\right)=\epsilon\left(K^2+1\right)
\left(\left|42\right\rangle+\epsilon\left|24\right\rangle\right), \\
&&S_5:\qquad
T^{(2)}\left(\left|14\right\rangle-\left|41\right\rangle\right)=-\left(\left(K^2-2K+1\right)+K\left(q^2+q^{-2}\right)
\right)\left(\left|14\right\rangle-\left|41\right\rangle\right)\\
&&\phantom{S_5:}\qquad
T^{(2)}\left(\left|23\right\rangle-\left|32\right\rangle\right)=-\left(K^2-2K+1\right)\left(\left|23\right\rangle-
\left|32\right\rangle\right)
\end{eqnarray}
Finally, denoting
\begin{equation}
\left|A\right\rangle=\left(\left|14\right\rangle+\left|41\right\rangle\right),\qquad
\left|B\right\rangle=\left(\left|23\right\rangle+\left|32\right\rangle\right)
\end{equation}
and setting
\begin{equation}
T^{(2)}\left(a\left|A\right\rangle+b\left|B\right\rangle\right)=v
\left(a\left|A\right\rangle+b\left|B\right\rangle\right),
\end{equation}
where
\begin{equation}
v=K^2+1+K\cdot f
\end{equation}
$f$ being $K$-independent (a function $f\left(q\right)$ of $q$
only), one obtains the constraints
\begin{equation}
a\left(\left(q+q^{-1}\right)^2-f\right)+b\cdot
2\left(q+q^{-1}\right)=0,\qquad a\cdot
2\left(q+q^{-1}\right)+b\left(2-f\right)=0.
\end{equation}
Hence
\begin{equation}
f=\frac 12\left(q^2+q^{-2}+4\right)\pm \frac
12\sqrt{\left(q+q^{-1}\right)^4+12\left(q+q^{-1}\right)^2+4}
\end{equation}
with corresponding $K$-independent values of $\left(a,b\right)$.
The sum of the eigenvalues given by (5.17)-(5.24) is
\begin{equation}
4\left(K+1\right)^2
\end{equation}
consistently with (5.8) for $r=2$. Our explicit results for $r=2$
not only shows how the basic properties (1), (2), (3) stated at
the beginning of this section are all realized systematically but
also how (3.14) is carried over, the subspaces now being paired
via
\begin{equation}
\left(1,2,3,4\right)\rightarrow \left(4,3,2,1\right)
\end{equation}

The chain Hamiltonian for any $N$ is given by (4.1)-(4.3) with
(1.4)-(1.7) giving $P_0'$ in $\dot{\hat{R}}_0=\dot{K}_0P_0'$. In
$\dot{K}_0$ of (4.5) now, from (1.3),
$e^\eta+e^{-\eta}=\left[N-1\right]+1$. For $r=2$, one obtains for
example
\begin{equation}
H^{(2)}=\dot{K}_0\left(P_0'+PP_0'P\right)=\dot{K}_0\left(\sum_{i,j=1}^N\left(q^{\rho_{i'}-\rho_j}+
q^{\rho_{i}-\rho_j'}\right)(ij)\otimes(i'j')\right).
\end{equation}
For $N=4$ this corresponds to
\begin{eqnarray}
&&\left(\dot{K}_0\right)^{-1}H^{(2)}=\left(q^{-2}+q^{2}\right)(11)\otimes
(44)+\left(q^{-1}+q\right)(12)\otimes (43)+
\left(q^{-1}+q\right)(13)\otimes (42)+2(14)\otimes (41)\nonumber\\
&&\phantom{\left(\dot{K}_0\right)^{-1}H^{(2)}=}\left(q^{-1}+q\right)(21)\otimes
(34)+2(22)\otimes (33)+2(23)\otimes (32)+\left(q^{-1}+q\right)(24)\otimes (31)\nonumber\\
&&\phantom{\left(\dot{K}_0\right)^{-1}H^{(2)}=}\left(q^{-1}+q\right)(31)\otimes
(24)+2(32)\otimes (23)+2(33)\otimes (22)+\left(q^{-1}+q\right)(34)\otimes (21)\\
&&\phantom{\left(\dot{K}_0\right)^{-1}H^{(2)}=}2(41)\otimes
(14)+\left(q^{-1}+q\right)(42)\otimes
(13)+\left(q^{-1}+q\right)(43)\otimes (12)+\left(q^{-2}+q^2
\right)(44)\otimes (11).\nonumber
\end{eqnarray}
Generalizations for $r>2$ can be written down systematically. If
the "spin" associated with the state $\left|i\right\rangle$ is
denoted as $\sigma_i$ then (4.1) along with structures analogous
to (5.27) implies transitions (if the states of two neighboring
sites have spins $\sigma_j$, $\sigma_{j'}$)
\begin{equation}
 \left(\sigma_j,\sigma_{j'}\right)\rightarrow  \left(\sigma_i,\sigma_{i'}\right)
\end{equation}
with evident $q$-dependent transition amplitudes corresponding to
the matrix elements of $H^{(r)}$ for order $r$. In particular if,
for example,
\begin{equation}
\left(\sigma_1,\sigma_2,\ldots,\sigma_{N-1},\sigma_N\right)=
\left(\frac{N-1}2,\frac{N-2}2,\ldots,-\frac{N-2}2,-\frac{N-1}2\right)
\end{equation}
then
\begin{equation}
\sigma_i+\sigma_{i'}=0\qquad \qquad (i'=N-i+1).
\end{equation}
For $N=3$ (as discussed in (4.26)-(4.28))
\begin{equation}
\left(\sigma_1,\sigma_2,\sigma_{3}\right)= \left(1,0,-1\right)
\end{equation}
and for $N=4$
\begin{equation}
\left(\sigma_1,\sigma_2,\sigma_{3},\sigma_4\right)=
\left(\frac{3}2,\frac{1}2,-\frac{1}2,-\frac{3}2\right)
\end{equation}
and so on.

If for two adjacent sites (including circular boundary
constraints) one has states $\left(\left|i\right\rangle,
\left|i'\right\rangle\right)$ they can be flipped to any pair
$\left(\left|j\right\rangle, \left|j'\right\rangle\right)$. Thus
such a flip can propagate along the chain for any state
$\left|j\right\rangle$ of the next site.

A thorough investigation of our class of models for arbitrary $N$
is beyond the scope of the present paper. We have however
indicated how the basic features studied for $N=3$ are carried
over as $N$ increases. Such properties are conserved due to the
specific structure of $P_0'$ as defined in (1.4)-(1.7).

We just mention finally that features parallel to those discussed
for $N=3$ in (3.43)-(3.49) can be carried over starting with the
multinomial expansion of
\begin{equation}
\left(x_1+x_2+x_3+\cdots+x_N\right)^r.
\end{equation}
Dimensions of invariant subspaces are obtained entirely
analogously.

\section{Potential for factorizable $S$-matrix $(N=3)$}
\setcounter{equation}{0}

As in Sec. 5 of Ref. 3 we construct the inverse Cayley transform
of the YB matrix which is also the $t^{(1)}\left(\theta\right)$
matrix (2.1) and given by (2.2) for $N=3$ for the class studied in
this paper. The role of this in providing the potential for
factorizable $S$-matrices can be found in various sources
\cite{R5,R6}. As explained and emphasized in Sec. 5 of Ref. 3 an
arbitrary normalization factor (denoted
$\lambda^{-1}\left(\theta\right)$) of $R\left(\theta\right)$ must
be introduced to start with for the inversion involved in the
transform to be well-defined. The explicit inversion in the first
factor of
\begin{equation}
-iV=\left(R\left(\theta\right)-\lambda\left(\theta\right)I\right)^{-1}\left(R\left(\theta\right)
+\lambda\left(\theta\right)I\right)
\end{equation}
will display admissible choices of $\lambda\left(\theta\right)$.
Defining
\begin{equation}
X\left(R\left(\theta\right)-\lambda\left(\theta\right)I\right)=I,\qquad
-iV=X\left(X^{-1}+2\lambda\left(\theta\right)I\right)=I+2\lambda\left(\theta\right)X,
\end{equation}
for $N=3$, (2.1) leads to (suppressing the argument $\theta$ in
notation below)
\begin{equation}
X\left|\begin{array}{ccccccccc}
   1-\lambda & 0 & 0 & 0 & 0& 0 & 0& 0 & 0\\
   0 & -\lambda & 0 & 1 & 0& 0 & 0& 0 & 0\\
   0 & 0 & K-\lambda & 0 & q^{1/2}K& 0 & 1+qK& 0 & 0\\
   0 & 1 & 0 & -\lambda & 0& 0 & 0& 0 & 0\\
   0 & 0 & q^{-1/2}K & 0 & 1+K-\lambda & 0 & q^{1/2}K& 0 & 0\\
   0 & 0 & 0 & 0 & 0& -\lambda & 0& 1 & 0\\
   0 & 0 & 1+q^{-1}K & 0 & q^{-1/2}K& 0 & K-\lambda& 0 & 0\\
   0 & 0 & 0 & 0 & 0& 1 & 0& -\lambda & 0\\
   0 & 0 & 0 & 0 & 0& 0 & 0& 0 & 1-\lambda\\
\end{array}\right|=I,
\end{equation}
Only the non-zero elements of $X$ will be given below. One obtains
easily
\begin{equation}
\left(X_{11},X_{99}\right)=\left(1-\lambda\right)^{-1},\qquad
\left(X_{22},X_{44},X_{66},X_{88}\right)=\frac\lambda{\left(1-\lambda^2\right)},
\qquad\left(X_{24},X_{42},X_{68},X_{86}\right)=\frac
1{\left(1-\lambda^2\right)}.
\end{equation}
These already show $\lambda \neq \pm 1$. For
$i=\left(3,5,7\right)$ one obtains the equations
\begin{equation}
\left(-X_{i3}\lambda+X_{i7}\right)+q^{-1/2}KZ_i=\delta_{i3},\qquad
X_{i5}\left(1-\lambda\right)+KZ_i=\delta_{i5},\qquad
\left(X_{i3}-\lambda X_{i7}\right)+Kq^{1/2}Z_i=\delta_{i7},
\end{equation}
where
\begin{eqnarray}
&&Z_i\equiv q^{1/2}X_{i3}+X_{i5}+q^{-1/2}X_{i7}.
\end{eqnarray}
The solutions for $i=\left(3,5,7\right)$ are respectively the
following ones. For $i=3$, $\left(X_{33},X_{35},X_{37}\right)$ are
given by
\begin{eqnarray}
&&X_{33}=\frac\lambda{1-\lambda^2}+q^{-1/2}\left(\frac{q+\lambda}{1+\lambda}\right)X_{35},\\
&&X_{37}=\frac 1{1-\lambda^2}+q^{-1/2}\left(\frac{1+q\lambda}{1+\lambda}\right)X_{35},\\
&&Z_{3}=\frac{q^{-1/2}\left(1+q\lambda\right)}{1-\lambda^2}+\left(\frac{q+q^{-1}+1+3\lambda}
{1+\lambda}\right)X_{35},\\
&&X_{35}\left(1-\lambda\right)+KZ_3=0.
\end{eqnarray}
The $K$-dependence is now explicit. The case $i=3,7$ are related
though the exchange of indices and inversion of $q$, namely
\begin{equation}
\left(3,7;q\right)\rightleftharpoons \left(7,3;q^{-1}\right)
\end{equation}
For $i=5$
\begin{eqnarray}
&&X_{53}\left(1+q\lambda\right)=X_{57}\left(q+\lambda\right),\\
&&X_{53}=-\frac {q^{-1/2}\left(q+\lambda\right)}{1-\lambda^2}+\frac{\left(q+\lambda\right)q^{-1/2}}{1+\lambda}X_{55},\\
&&X_{57}=-\frac {q^{-1/2}\left(1+q\lambda\right)}{1-\lambda^2}+\frac{\left(1+q\lambda\right)q^{-1/2}}{1+\lambda}X_{55},\\
&&Z_{5}=-\frac
{q+q^{-1}+2\lambda}{1-\lambda^2}+\frac{3\lambda+1+q+q^{-1}}
{1+\lambda}X_{55},\\
&&X_{55}\left(1-\lambda\right)+KZ_5=1.
\end{eqnarray}
Now (6.19), (6.20) gives directly $X_{55}$. Next (6.17), (6.18)
give $X_{53}$, $X_{57}$. Finally, we obtain
\begin{equation}
X=\left|\begin{array}{ccccccccc}
   \frac 1{1-\lambda} & 0 & 0 & 0 & 0& 0 & 0& 0 & 0\\
   0 & \frac \lambda{1-\lambda^2} & 0 & \frac 1{1-\lambda^2} & 0& 0 & 0& 0 & 0\\
   0 & 0 & {\sf A} & 0 & {\sf B}& 0 & {\sf C}& 0 & 0\\
   0 & \frac 1{1-\lambda^2} & 0 & \frac \lambda{1-\lambda^2} & 0& 0 & 0& 0 & 0\\
   0 & 0 & {\sf D} & 0 & {\sf E}& 0 & {\sf B}& 0 & 0\\
   0 & 0 & 0 & 0 & 0& \frac \lambda{1-\lambda^2} & 0& \frac 1{1-\lambda^2} & 0\\
   0 & 0 & {\sf F} & 0 & {\sf D}& 0 & {\sf A}& 0 & 0\\
   0 & 0 & 0 & 0 & 0& \frac 1{1-\lambda^2} & 0& \frac \lambda{1-\lambda^2} & 0\\
   0 & 0 & 0 & 0 & 0& 0 & 0& 0 & \frac 1{1-\lambda}\\
\end{array}\right|,
\end{equation}
where
\begin{eqnarray}
&&{\sf
A}=\frac{\left(\lambda^2-\lambda-2K\lambda+K\right)q}{\left(1-\lambda\right)\left(q\lambda^2-3qK\lambda-q^2K-K-q
K-q\right)}=X_{33}=X_{77},\nonumber\\
&&{\sf
B}=\frac{\left(q\lambda+1\right)q^{1/2}K}{\left(1-\lambda\right)\left(q\lambda^2-3qK\lambda-q^2K-K-q
K-q\right)}=X_{35}=X_{57},\nonumber\\
&&{\sf
C}=\frac{\left(\lambda+qK\lambda-qK-1-K\right)q}{\left(1-\lambda\right)\left(q\lambda^2-3qK\lambda-q^2K-K-q
K-q\right)}=X_{37},\nonumber\\
&&{\sf
D}=\frac{K\left(\lambda+q\right)q^{1/2}}{\left(1-\lambda\right)\left(q\lambda^2-3qK\lambda-q^2K-K-q
K-q\right)}=X_{53}=X_{75},\nonumber\\
&&{\sf
E}=\frac{q\lambda^2-2Kq\lambda-q^2K-K-q}{\left(1-\lambda\right)\left(q\lambda^2-3qK\lambda-q^2K-K-q
K-q\right)}=X_{55},\nonumber\\
&&{\sf
F}=\frac{-q-qK+q\lambda-K+K\lambda}{\left(1-\lambda\right)\left(q\lambda^2-3qK\lambda-q^2K-K-q
K-q\right)}=X_{73},
\end{eqnarray}
where $\lambda\neq \pm 1, \frac 12\left[3K\pm
\sqrt{9K^2+4\left(q+1+q^{-1}\right)K+4}\right]$.

From $X$, $V$ is obtained as indicated in (6.2). Expressing it as
\begin{equation}
V=\sum_{ab,cd}V_{(ab,cd)}\left(ab\right)\otimes \left(cd\right)
\end{equation}
The corresponding fermionic Lagrangian should be
\begin{equation}
{\cal L}=\int
dx\left(i\overline{\psi}_a\gamma_\nu\partial_\nu\psi_a-g\left(\overline{\psi}_a\gamma_\nu\psi_c\right)
V_{ab,cd}\left(\overline{\psi}_b\gamma_\nu\psi_d\right)\right),
\end{equation}
The scalar Lagrangian can be obtained analogously. Such
Lagrangians correspond to $S$-matrices factorizable into two
particles scattering independently of the chosen order of the
latter ones.

\section{Discussion}
\setcounter{equation}{0}

In Ref. 3 and in the present paper we have studied two different
classes of statistical models. Certain aspects of the respective
transfer matrices are strikingly contrasted. Such a major
difference is in the number of parameters. The first model is
indeed multiparameter. One has $\frac
12\left(N+3\right)\left(N-1\right)$ free parameters
$\left(N=3,4,\ldots\right)$. Here the only parameter is $q$
appearing in the braid matrix given by (1.2)-(1.7) and in
$K\left(\theta\right)$ as defined by (1.17)-(1.20). The structures
of the eigenvalues of the respective transfer matrices are also
quite different. In Ref. 3 we obtained single exponentials as
eigenvalues, the exponent being a sum of the free parameters
multiplied by $\theta$. Here we have $r$-th order polynomials in
$K\left(\theta\right)$ for the eigenvalues of
$T^{(r)}\left(\theta\right)$. There are other differences. But
analogies and common features are also remarkable:

\begin{description}
    \item[(a)] In both case $\hbox{Tr}\left(T^{(r)}\left(\theta\right)\right)$ is
    obtained quite simply for all $r$ (though the structures are
    different). In (6.1) of Ref. 3 we obtained (for $N=2p-1$)
\begin{equation}
\hbox{Tr}\left(T^{(r)}\left(\theta\right)\right)=2\left(e^{rm_{11}^{(+)}\theta}+e^{rm_{22}^{(+)}\theta}+\cdots+e^{rm_{p-1,p-1}^{(+)}
\theta}\right)+1
\end{equation}
the $m_{ii}^{(+)}$ being a subset of the free parameters. Here
(for $N=3,4,\ldots$) the corresponding result (5.9) is
\begin{equation}
\hbox{Tr}\left(T^{(r)}\left(\theta\right)\right)=N\left(K\left(\theta\right)+1\right)^r,
\end{equation}
where $K\left(\theta\right)$ is given by (1.3) and (1.17).
    \item[(b)] In both cases the $N^r$ dimensional base space of
    $T^{(r)}\left(\theta\right)$ breaks up into closed subspaces
    of lower dimensions. The definitions of these subspaces have
    some differences however. The relevant definitions in Ref. 3
    should be compared to (3.2)-(3.6) here and their
    generalization in Sec. 5.
    \item[(c)] In each subspace $S_n$ the circular permutation of
    states labels as formulated in (3.9)-(3.19) leads to a
    further reduction of dimension in constructing eigenstates by
    splitting $S_n$ again into subsets corresponding to the eigenstates
    of the operator (CP) of circular permutations. This involves
    a crucial role of the roots of unity, ($\omega^r=1$ for
    $T^{(r)}\left(\theta\right)$) in the construction of
    eigenstates. The role of roots of unity was also crucial in
    Ref. 3 though they were implemented in a slightly different
    fashion (corresponding to the difference in labeling states).
\end{description}

In both cases the "two-step reduction" (via (b) and (c)) in the
effective dimension of the basis in construction of eigenstates
has been emphasized (see the formulation of Sec. 3). The
exponential increase in dimension with $r$ $\left(e^{\left(\ln
N\right)r}\right)$ is replaced in actual construction by a
relatively moderate polynomial one. Thus for $N=3$ and $r=4$ we
have to solve here at most a set of 5 simultaneous linear
equations (App. A) though now $N^r$ is $3^4=81$. This reduction of
the problem to a  relatively low number of linear equations should
be contrasted to the implementation of algebraic Bethe ansatz
\cite{R6,R7,R8}. For the latter one has to solve complex nonlinear
equations whose number increases along with $N$.

In the preceding sections (particularly in Sec. 3 for $N=3$ and in
sec. 5 for $N>3$), we have formulated carefully the crucial
properties, basic features of models corresponding to the braid
matrices presented in Ref. 2. Exploiting such properties we have
constructed eigenstates and eigenvalues of
$T^{(r)}\left(\theta\right)$ for $N=3$, $r=\left(1,2,3,4\right)$
(App. A). Certain related features for all $r$ have been
formulated at the end of Sec. 3. Chain Hamiltonians and potentials
for factorizable $S$-matrices have been studied (sections 4 and
6).

Further explorations in several directions are evidently
desirable. Reflection equations \cite{R9,R10} and correlation
functions \cite{R11,R12} should be studied. More basically one may
try to elucidate the relevance of the star-triangle relations
\cite{R4} encoded in our class of braid matrices to specific
contexts. We hope to undertake such studies elsewhere.

  \vskip 0.5cm

\noindent{\bf Acknowledgments:} {\em One of us (BA) wants to thank
Patrick Mora, Pierre Collet and Paul Sorba for precious help. He
is also very grateful to the members of the CPHT group of Ecole
Polytechnique for their warm hospitality. This work is supported
by a grant of CMEP program under number 04MDU615.}

\begin{appendix}

\section{\LARGE Eigenstates and Eigenvalues of $T^{(r)}\left(\theta\right)$ for
$r=1,2,3,4$ $\left(N=3\right)$} \setcounter{equation}{0}

We start by noting the following points:

\begin{description}
    \item[(1)] For each case the subscript $n$ of $S_n$ denotes
    the sum of the state labels (see discussion from (3.3) to
    (3.6)).
    \item[(2)] For each $r$ we present results only upto $S_{2r}$. The remaining
    subspaces $\left(S_{2r+1},\ldots,S_{3r}\right)$ are then
    obtained implementing (3.14), (3.15), (3.16).
    \item[(3)] For different subspaces we often repeat the some
    notations for states. Since $T^{(r)}$ does not couple such
    spaces no confusion is likely.
 \item[(4)] The notation $K$ and $\omega$ correspond to (1.17) and (3.11) respectively.
 (CP) denotes circular permutations.
\end{description}

 \subsection{$r=1$} The three states directly furnish the spectrum, each
 being a 1-dimensional subspace
\begin{eqnarray}
&&T^{(1)}\left(\theta\right)\left(\left|1\right\rangle,\left|2\right\rangle,\left|3\right\rangle\right)
=(1+K)\left(\left|1\right\rangle,\left|2\right\rangle,\left|3\right\rangle\right)\\
&&\hbox{Tr}\left(T^{(1)}\left(\theta\right)\right)=3(1+K)
\end{eqnarray}

 \subsection{$r=2$} The (CP) eigenstates constructed as in (3.10) to (3.13), with $\omega^2=1$,
give
\begin{eqnarray}
&&S_2:\qquad A_1=\left|11\right\rangle\nonumber\\
&&S_3:\qquad A_{\pm 1}=\left|12\right\rangle\pm\left|21\right\rangle \nonumber\\
&&S_4:\qquad A_1=\left|22\right\rangle,\qquad\qquad B_{\pm
1}=\left|13\right\rangle\pm\left|31\right\rangle
\end{eqnarray}
One obtains ($K$ being $K\left(\theta\right)$)
\begin{eqnarray}
&&S_2:\,\, T^{(2)}\left(\theta\right)A_1=\left(K^2+1\right)A_1\\
&&S_3:\,\, T^{(2)}\left(\theta\right)A_{\pm 1}=\pm\left(K^2+1\right)A_{\pm 1}\\
&&S_4:\,\, T^{(2)}\left(\theta\right)B_{-1}=-\left(K^2+\left(q+q^{-1}-2\right)K+1\right)B_{-1}\nonumber\\
&&\phantom{S_4:}\,\,
T^{(2)}\left(\theta\right)\left(B_{1}-\left(q^{1/2}+q^{-1/2}\right)A_1\right)=
\left(K^2+1\right)\left(B_{1}-\left(q^{1/2}+q^{-1/2}\right)A_1\right)\\
&&\phantom{S_4:}\,\,
T^{(2)}\left(\theta\right)\left(\left(q^{1/2}+q^{-1/2}\right)B_{1}+2A_1\right)=
\left(K^2+\left(q+q^{-1}+4\right)K+1\right)\left(\left(q^{1/2}+q^{-1/2}\right)B_{1}+2A_1\right).\nonumber
\end{eqnarray}
Also
\begin{equation}
\left(S_5,S_6\right)\rightleftharpoons \left(S_3,S_2\right)
\end{equation}
according to (3.14)-(3.16). Summing over all subspaces
$\left(S_2,\ldots, S_6\right)$,
\begin{equation}
\hbox{Tr}\left(T^{(1)}\left(\theta\right)\right)=3\left(1+K\right)^2
\end{equation}
consistently with (2.9). We have not uniformly normalized the
states. Thus $\left\langle A_1|A_1\right\rangle=1$ and
$\left\langle B_{\pm 1}|B_{\pm 1}\right\rangle=2$. This is crucial
to the orthogonality
\begin{equation}
\left\langle
B_1-\left(q^{1/2}+q^{-1/2}\right)A_1|\left(q^{1/2}+q^{-1/2}\right)B_1+2A_1\right\rangle=0.
\end{equation}
This point displayed here for this simple case will not be
repeated in cases to follow.

 \subsection{$r=3$} Here
\begin{equation}
\omega=\left(1,e^{i\frac{2\pi}3},e^{i\frac{2\pi}3\cdot 2}\right)
\end{equation}
and (CP) eigenstates for $\left(S_3,S_4,S_5,S_6\right)$ are
\begin{eqnarray}
&&S_3:\,\, A_1=\left|111\right\rangle\\
&&S_4:\,\, A_{\omega}=\left|112\right\rangle+\omega\left|211\right\rangle+\omega^2\left|121\right\rangle\\
&&S_5:\,\,
A_{\omega}=\left|113\right\rangle+\omega\left|311\right\rangle+\omega^2\left|131\right\rangle,\,\,
B_{\omega}=\left|122\right\rangle+\omega\left|212\right\rangle+\omega^2\left|221\right\rangle\\
&&S_6:\,\, A_{1}=\left|222\right\rangle,\,\,
B_{\omega}=\left|123\right\rangle+\omega\left|312\right\rangle+\omega^2\left|231\right\rangle,
\nonumber\\
&&\phantom{S_6:\,\, }
C_{\omega}=\left|321\right\rangle+\omega\left|132\right\rangle+\omega^2\left|213\right\rangle
=\left(1\rightleftharpoons 3\right)B_\omega.
\end{eqnarray}
Also
\begin{equation}
\left(1\rightleftharpoons
3\right)\left(S_3,S_4,S_5\right)=\left(S_9,S_8,S_7\right)
\end{equation}

The $T^{(3)}\left(\theta\right)$ eigenstates are now obtained as
follows:
\begin{eqnarray}
&&S_3:\qquad T^{(3)}\left(\theta\right)A_1=\left(K^3+1\right)A_1\\
&&S_4:\qquad T^{(3)}\left(\theta\right)A_\omega=\left(K^3\omega+\omega^2\right)A_\omega\\
&&S_5:\qquad
T^{(3)}\left(\theta\right)\left(aA_\omega+bB_\omega\right)=v\left(aA_\omega+bB_\omega\right)
\end{eqnarray}
{\bf Solutions}:
\begin{eqnarray}
&&\hbox{\bf (1)}\,\,\,\, \left(a,b\right)=\left(q^{1/2}+\omega q^{-1/2},1\right),\nonumber\\
&&\phantom{\hbox{\bf
(1)}\,\,\,\,}v=\omega^2K^3+\left(q+q^{-1}\right)\left(\omega^2K^2+\omega
K\right)+\left(1+\omega+\omega^2\right)\left(K^2+K\right)+\omega\\
&&\hbox{\bf (2)}\,\,\,\,
\left(a,b\right)=\left(1,-\left(q^{1/2}+\omega^2
q^{-1/2}\right)\right),\qquad v=K^3\omega^2+\omega.
\end{eqnarray}

\paragraph{$\diamond$ \underline{$S_6$}:} For $\omega=e^{\pm i\frac{2\pi}3}$, $A_1$ is
decoupled. Set
\begin{equation}
T^{(3)}\left(\theta\right)\left(bB_\omega+cC_\omega\right)=v\left(bB_\omega+cC_\omega\right).
\end{equation}
{\bf Solutions}:
\begin{eqnarray}
&&\hbox{\bf (1)}\,\,\,\, \left(b,c\right)=\left(q,-1\right),\qquad v=K^3\omega^2+\omega\\
&&\hbox{\bf (2)}\,\,\,\,\left(b,c\right)=\left(1,q\right),\qquad
v=K^3\omega^2+\omega+\left(q+q^{-1}\right)\left(K^2\omega^2+K\omega\right).
\end{eqnarray}
For the values of $\omega$ $\left(\neq 1\right)$, with
$\omega+\omega^2=-1$ the sum of eigenvalues
\begin{equation}
\sum
v=-2\left(K^3+1\right)-\left(q+q^{-1}\right)\left(K^2+K\right).
\end{equation}
For $\omega=1$, $T^{(3)}$ couples $\left(A_1,B_1,C_1\right)$. Set
\begin{equation}
T^{(3)}\left(\alpha A_1+\beta B_1+\gamma
C_1\right)=v_1\left(\alpha A_1+\beta B_1+\gamma C_1\right).
\end{equation}
{\bf Solutions}:
\begin{eqnarray}
&&\hbox{\bf (1)}\,\,\,\, \left(\alpha,\beta,\gamma\right)=\left(0,q^{1/2},-q^{-1/2}\right),\qquad
v_1=\left(K+1\right)\left(K^2-K+1\right)\\
&&\hbox{\bf (2)}\,\,\,\, \left(\alpha,\beta,\gamma\right)=\left(-\left(q+q^{-1}\right),q^{-1/2},q^{1/2}\right),\qquad
v_1=\left(K+1\right)\left(K^2-K+1\right)\\
&&\hbox{\bf (3)}\,\,\,\,
\left(\alpha,\beta,\gamma\right)=\left(3,q^{-1/2},q^{1/2}\right),\qquad
v_1=\left(K+1\right)\left(\left(K^2-K+1\right)+K\left(q+q^{-1}+3\right)\right).
\end{eqnarray}
Concerning orthogonality note that
\begin{equation}
\left\langle A_1|A_1\right\rangle=1,\qquad \left\langle
B_1|B_1\right\rangle=\left\langle C_1|C_1\right\rangle=3.
\end{equation}
The sum of the eigenvalues over $S_6$ is
\begin{equation}
\left(\sum v\right)_{S_6}=K^3+1+3K\left(K+1\right).
\end{equation}
The results for $\left(S_7,S_8,S_9\right)$ are obtained, as usual,
directly from those of $\left(S_5,S_4,S_3\right)$ respectively.

Summing over all the subspaces $\left(S_3,\ldots,S_9\right)$ one
obtains (all explicit $q$-dependence canceling consistently with
(2.9))
\begin{eqnarray}
&&Tr\left(T^{(3)}\left(\theta\right)\right)=\left(K^3+1\right)+\left(K^3+1\right)+
3\left(K^2+K\right)+3\left(K^2+K\right)+\left(K^3+1\right)+3\left(K^2+K\right)\nonumber\\
&&\phantom{Tr\left(T^{(3)}\left(\theta\right)\right)}=3\left(K+1\right)^3.
\end{eqnarray}

 \subsection{$r=4$} Here
\begin{equation}
\omega=\left(1,e^{i\frac{2\pi}4},e^{i\frac{2\pi}4\cdot
2},e^{i\frac{2\pi}4\cdot 3}\right)=\left(1,i,-1,-i\right).
\end{equation}
Of the invariant subspaces we consider $\left(S_4,S_5,S_6,
S_7,S_8\right)$. One obtains the results of the remaining ones via
$\left(1,2,3;q\right)\leftrightarrow \left(3,2,1;q^{-1}\right)$ as
\begin{equation}
\left(S_9,S_{10},S_{11},S_{12}\right)\rightleftharpoons
\left(S_7,S_6,S_5, S_4\right).
\end{equation}
For brevity and simplicity, we will recapitulate our results in
the following tables:

\smallskip\smallskip\smallskip\smallskip\smallskip

\centerline{{\bf Table 1}: {\it CP eigenstates for $r=4$}}
$$\begin{array}{llc}
\hline&\nonumber\\
\hbox{Subspace}\qquad &\hbox{CP Eigenstates}& \hbox{Dimension}\qquad\nonumber \\
\hline &&\nonumber\\
S_4&  \left|1111\right\rangle & 1   \nonumber\\
&&\nonumber\\
S_5
&A_\omega=\left|1112\right\rangle+\omega\left|2111\right\rangle+\omega^2\left|1211\right\rangle+\omega^3
\left|1121\right\rangle & 4 \nonumber\\
&&\nonumber\\
S_6 & A_{\pm
1}=\left|1212\right\rangle\pm\left|2121\right\rangle & 10\nonumber\\
&B_{\omega}=\left|1113\right\rangle+\omega\left|3111\right\rangle+\omega^2\left|1311\right\rangle
+\omega^3\left|1131\right\rangle&\nonumber\\
&C_{\omega}=\left|1122\right\rangle+\omega\left|2112\right\rangle+\omega^2\left|2211\right\rangle
+\omega^3\left|1221\right\rangle&\nonumber\\
&\nonumber\\
S_7&A_{\omega}=\left|1222\right\rangle+\omega\left|2122\right\rangle+\omega^2\left|2212\right\rangle
+\omega^3\left|2221\right\rangle & 16\nonumber\\
&B_{\omega}=\left|1123\right\rangle+\omega\left|3112\right\rangle+\omega^2\left|2311\right\rangle
+\omega^3\left|1231\right\rangle &\nonumber\\
&C_{\omega}=\left|1132\right\rangle+\omega\left|2113\right\rangle+\omega^2\left|3211\right\rangle
+\omega^3\left|1321\right\rangle &\nonumber\\
&D_{\omega}=\left|1213\right\rangle+\omega\left|3121\right\rangle+\omega^2\left|1312\right\rangle
+\omega^3\left|2131\right\rangle&\nonumber\\
&&\nonumber\\
S_8 & A_1=\left|2222\right\rangle & 19\nonumber\\
&B_{\pm 1}=\left|1313\right\rangle\pm \omega\left|3131\right\rangle&\nonumber\\
&C_{\omega}=\left|1133\right\rangle+\omega\left|3113\right\rangle+\omega^2\left|3311\right\rangle
+\omega^3\left|1331\right\rangle &\nonumber\\
&D_{\omega}=\left|1223\right\rangle+\omega\left|3122\right\rangle+\omega^2\left|2312\right\rangle
+\omega^3\left|2231\right\rangle &\nonumber\\
&E_{\omega}=\left|3221\right\rangle+\omega\left|1322\right\rangle+\omega^2\left|2132\right\rangle
+\omega^3\left|2213\right\rangle &\nonumber\\
&F_{\omega}=\left|1232\right\rangle+\omega\left|2123\right\rangle+\omega^2\left|3212\right\rangle
+\omega^3\left|2321\right\rangle&\nonumber\\
\hline
\end{array}
$$

\newpage

\centerline{{\bf Table 2}: {\it Eigenstates and eigenvalues for
$r=4$}}
$$\begin{array}{lll}
\hline&\nonumber\\
 &\hbox{Eigenvalues}&\hbox{Eigenstates}\nonumber \\
\hline &\nonumber\\
S_4:& K^4+1  &  \left|1111\right\rangle\nonumber\\
&\nonumber\\&\nonumber\\
S_5: & \omega^3K^4+\omega & A_\omega\nonumber\\
&\nonumber\\&\nonumber\\
S_6: & \pm\left(K^4+1\right) & A_{\pm 1}\nonumber\\
&&\nonumber\\
& \omega^3K^4+\omega & B_\omega-\frac{q+\omega^3}{\sqrt{q}}C_\omega\nonumber\\
&&\nonumber\\
&\omega^3\left[K^4+(q+1+\omega+\omega^3+q^{-1})K^3+\right.&\nonumber\\
&\left.\omega^3(q+1+\omega+\omega^3+q^{-1})K^2+\right. &
B_\omega+\frac{\sqrt{q}}{q+\omega}C_\omega\nonumber\\
&\left.\omega^2(q+1+\omega+\omega^3+q^{-1})K+\omega^2\right]&\nonumber\\
&&\nonumber\\&&\nonumber\\
S_7: & \omega^3K^4+\omega & A_\omega-\omega\sqrt{q}B_\omega-\frac{1}{\sqrt{q}}C_\omega,\nonumber\\
& &
-\left(q^3+\omega^2\right)A_\omega+q\sqrt{q}\left(q+1+q^{-1}\right)D_\omega+\nonumber\\
&&\omega\sqrt{q}\left(\omega^2-q^2-q\right)B_\omega+\sqrt{q}\left(q^2-\omega^2q-\omega^2\right)C_\omega\nonumber\\
&&\nonumber\\
&\omega^3K^4+\omega+\left(\omega^3K^3+\omega
K\right)\left(q+1+\omega+\omega^3+q^{-1}\right)
&\frac{1+\omega}{\sqrt{q}}A_{\omega}+\frac
{\omega^2}{q}B_{\omega}+C_{\omega}+\omega\left(\frac
{q+\omega}q\right)D_{\omega}\nonumber\\
&+K^2\omega^2\left(q+1+\omega+\omega^3+q^{-1}\right)  &\nonumber\\
&&\nonumber\\
&\omega^3K^4+\omega+\left(\omega^3K^3+\omega
K\right)\left(q+1-\omega-\omega^3+q^{-1}\right)
&\frac{1-\omega}{\sqrt{q}}A_{\omega}-\frac
{\omega^2}{q}B_{\omega}+C_{\omega}-\omega\left(\frac
{q-\omega}q\right)
D_{\omega}\nonumber\\
&+K^2\omega^2\left(q+1-\omega-\omega^3+q^{-1}\right) &\nonumber\\
&&\nonumber\\&&\nonumber\\
S_8: & K^4+1 &
F_1,\; 2A_1+C_1-\sqrt{q}D_1-\frac 1{\sqrt{q}}E_1,\nonumber\\
&&
2\sqrt{q}\left(q^2+q-2+q^{-1}+q^{-2}\right)A_1+\nonumber\\
&&2\sqrt{q}\left(q+2+q^{-1}\right)B_1-\sqrt{q}\left(q+4+q^{-1}\right)C_1+\nonumber\\
&&\left(q^2+2q-3-2q^{-1}-3\right)D_1\nonumber\\
&&-\left(2q^2+3q-2-q^{-1}\right)E_1\nonumber\\
&& \nonumber\\
&-K^4-1& F_{-1},\nonumber\\
&&
B_{-1}+\left(\frac{q^2-1}{2q}\right)C_{-1}-\left(\frac{q^2+1}{2\sqrt{q}}\right)D_{-1}+
\left(\frac{q^2+1}{2q\sqrt{q}}\right)E_{-1}\nonumber\\
&& \nonumber\\
&\mp iK^4\pm i& F_{\pm i},\; \sqrt{q}C_{\pm i}-D_{\pm i}+E_{\pm i}\nonumber\\
&& \nonumber\\
&\mp i\left[K^4+\left(q+1+q^{-1}\right)K^3\pm 3iK^2-\right.& -\left(\frac{q+1}{\sqrt{q}}\right)C_{\pm i}-D_{\pm i}+
E_{\pm i}\nonumber\\
&\left.\left(q+1+q^{-1}\right)K-1\right]& \nonumber\\
&\nonumber\\
&\mp i\left[K^4+\left(q+1+q^{-1}\right)K^3\mp
i\left(q+1+q^{-1}\right)K^2-\right.& C_{\pm i}+\frac
{1+2q}{\left(q-1\right)\sqrt{q}}D_{\pm i}
+\frac{(q+2)\sqrt{q}}{q-1}E_{\pm i}\nonumber\\
&\left.\left(q+1+q^{-1}\right)K-1\right]& \nonumber\\
\end{array}
$$

$$\begin{array}{lll}
&K^4+\left(q+3+q^{-1}\right)K^3+\left(q+3+q^{-1}\right)K^2+&
-4A_1+2B_1+C_1-\frac 1{\left(q+1\right)\sqrt{q}}D_1
-\frac{q\sqrt{q}}{q+1}E_1\nonumber\\
&\left(q+3+q^{-1}\right)K+1& \nonumber\\
&\nonumber\\
&-\left[K^4+\left(q-1+q^{-1}\right)K^3-\left(q-1+q^{-1}\right)K^2+\right.&
-\frac{2\left(q-1\right)}{q+1}B_{-1}+C_{-1}-\frac
1{\left(q+1\right)\sqrt{q}}D_{-1}
-\frac{q\sqrt{q}}{q+1}E_{-1}\nonumber\\
&\left.\left(q-1+q^{-1}\right)K+1\right]& \nonumber\\
\end{array}
$$
$\star$ There exist also four others eigenvectors $a
A_1+bB_1+cC_1+dD_1+eE_1$, $\alpha A_1+\beta B_1+\gamma C_1+\delta
D_1+\eta E_1$, $b'B_{-1}+c'C_{-1}+d'D_{-1}+e'E_{-1}$, $\beta'
B_{-1}+\gamma' C_{-1}+\delta' D_{-1}+\eta' E_{-1}$ associated
respectively to the eigenvalues $v_1$, $v_2$, $v_1'$ and $v_2'$,
which have complicated forms (these results have been obtained by
using a Maple program):
\begin{eqnarray}
&&v_1=\frac
12\left(2K^{4}+3K^3\left(q+1+q^{-1}\right)+K^2\left({q}^2+2{q}+13+2q^{-1}+q^{-2}\right)
+3K\left(q+1+q^{-1}\right)+2+\right.\nonumber\\
&&\phantom{v_1=}\left.K\left(K^2+\left(q+1+q^{-1}\right)K+1\right)\sqrt{q^2+2q+43+2q^{-1}+q^{-2}}\right),\nonumber\\
&&v_2=\frac
12\left(2K^{4}+3K^3\left(q+1+q^{-1}\right)+K^2\left({q}^2+2{q}+13+2q^{-1}+q^{-2}\right)
+3K\left(q+1+q^{-1}\right)+2-\right.\nonumber\\
&&\phantom{v_2=}\left.K\left(K^2+\left(q+1+q^{-1}\right)K+1\right)\sqrt{q^2+2q+43+2q^{-1}+q^{-2}}\right),\nonumber\\
&&v_1'=\frac
12\left(-2K^{4}-3K^3\left(q+1+q^{-1}\right)-K^2\left({q}^2+2{q}+1+2q^{-1}+q^{-2}\right)
-3K\left(q+1+q^{-1}\right)-2+\right.\nonumber\\
&&\phantom{v_1'=}\left.K\left(K^2+\left(q+1+q^{-1}\right)K+1\right)\sqrt{q^2+2q-5+2q^{-1}+q^{-2}}\right),\nonumber\\
&&v_2'=\frac
12\left(-2K^{4}-3K^3\left(q+1+q^{-1}\right)-K^2\left({q}^2+2{q}+1+2q^{-1}+q^{-2}\right)
-3K\left(q+1+q^{-1}\right)-2-\right.\nonumber\\
&&\phantom{v_2'=}\left.K\left(K^2+\left(q+1+q^{-1}\right)K+1\right)\sqrt{q^2+2q-5+2q^{-1}+q^{-2}}\right).\nonumber\\
&&
\end{eqnarray}
$\underline{\phantom{XXXXXXXXXXXXXXXXXXXXXXXXXXXXXXXXXXXXXXXXXXXXXXXXX}}$

\smallskip\smallskip\smallskip\smallskip\smallskip

The sums of the eigenvalues over $S_4$, $S_5$, $S_6$, $S_7$ and
$S_8$ are respectively
\begin{eqnarray}
&&\sum_{S_4}v=K^4+1,\nonumber\\
&&\sum_{S_5}v=0,\nonumber\\
&&\sum_{S_6}v=4K^3+4K,\nonumber\\
&&\sum_{S_7}v=0,\nonumber\\
&&\sum_{S_8}v=K^4+4K^3+18K^2+4K+1.
\end{eqnarray}
The results for $\left(S_9,S_{10},S_{11},S_{12}\right)$ are
obtained directly from those of $\left(S_7,S_6,S_5,S_4\right)$
respectively. Summing over all subspaces
$\left(S_4,\ldots,S_{12}\right)$ one obtains
\begin{equation}
\hbox{Tr}\left(T^{(4)}\left(\theta\right)\right)=2\times\left(K^4+1\right)+
2\times \left(4K^3+4K\right)+
1\times\left(K^4+4K^3+18K^2+4K+1\right)=3\left(K+1\right)^4.
\end{equation}

\section{\LARGE $\hat{R}tt$-Algebra} \setcounter{equation}{0}

We present below, for $N=3$, the constraints on the blocks
$t_{ij}\left(\theta\right)$ of the transfer matrix following from
\begin{equation}
\hat{R}\left(\theta-\theta'\right)t\left(\theta\right)\otimes
t\left(\theta'\right)=t\left(\theta'\right)\otimes\left(\theta\right)\hat{R}\left(\theta-\theta'\right).
\end{equation}
We use below the notations
\begin{equation}
\left(t\left(\theta\right),t\left(\theta'\right),K\left(\theta-\theta'\right)\right)\equiv
\left(t,t',K''\right)
\end{equation}
In terms of $P_0'$ defined by (1.4) with $\left(i,j\right)$ and
$\left(\rho_i,\rho_j\right)$ corresponding to $N=3$ and
$K\left(\theta\right)$ of (1.17)-(1.18), (B.1) now is (maintaining
the notation $P_0'$ unrelated to $\theta$, $\theta'$)
\begin{equation}
\left(I+K''P_0'\right)\left(t\otimes t'\right)=\left(t'\otimes
t\right)\left(I+K''P_0'\right),
\end{equation}
where
\begin{eqnarray}
&&P_0'=q^{-1}\left(11\right)\otimes
\left(33\right)+q^{-1/2}\left(12\right)\otimes
\left(32\right)+\left(13\right)\otimes
\left(31\right)+q^{-1/2}\left(21\right)\otimes
\left(23\right)+\left(22\right)\otimes
\left(22\right)\nonumber\\
&&\phantom{P_0'=}+q^{1/2}\left(23\right)\otimes
\left(21\right)+\left(31\right)\otimes
\left(13\right)+q^{1/2}\left(32\right)\otimes
\left(12\right)+q\left(33\right)\otimes \left(11\right)
\end{eqnarray}
This leads to a set of 36 relations independent of $K''$, namely
\begin{equation}
t_{ij}t'_{kl}=t'_{ij}t_{kl}
\end{equation}
where for $(ij)=(11),(12),(13)$ respectively
\begin{equation}
(kl)=(11,12,21,22),\,(11,13,21,23),\,(12,13,22,23)
\end{equation}
and similarly for $(ij)=(21),(22),(23)$
\begin{equation}
(kl)=(11,12,31,32),\,(11,13,31,33),\,(12,13,32,33)
\end{equation}
and for $(ij)=(31),(32),(33)$
\begin{equation}
(kl)=(21,22,31,32),\,(21,23,31,33),\,(22,23,32,33)
\end{equation}

To present the $K''$ dependent constraints we first define
\begin{eqnarray}
&&X_1=q^{-1/2}t_{11}t'_{31}+t_{21}t'_{21}+q^{1/2}t_{31}t'_{11},\qquad
X_2=q^{-1/2}t_{11}t'_{32}+t_{21}t'_{22}+q^{1/2}t_{31}t'_{12},\nonumber\\
&&X_3=q^{-1/2}t_{11}t'_{33}+t_{21}t'_{23}+q^{1/2}t_{31}t'_{13},\qquad
X_4=q^{-1/2}t_{12}t'_{31}+t_{22}t'_{21}+q^{1/2}t_{32}t'_{11},\nonumber\\
&&X_5=q^{-1/2}t_{12}t'_{32}+t_{22}t'_{22}+q^{1/2}t_{32}t'_{12},\qquad
X_6=q^{-1/2}t_{12}t'_{33}+t_{22}t'_{23}+q^{1/2}t_{32}t'_{13},\nonumber\\
&&X_7=q^{-1/2}t_{13}t'_{31}+t_{23}t'_{21}+q^{1/2}t_{33}t'_{11},\qquad
X_8=q^{-1/2}t_{13}t'_{32}+t_{23}t'_{22}+q^{1/2}t_{33}t'_{12},\nonumber\\
&&X_9=q^{-1/2}t_{13}t'_{33}+t_{23}t'_{23}+q^{1/2}t_{33}t'_{13}
\end{eqnarray}
and a set
\begin{equation}
\left(Y_1,Y_2,\ldots,Y_9\right)
\end{equation}
which is obtained by transposing the indices of each term on the
right of (B.9) and also the order of
$\left(\theta,\theta'\right)$. Thus
\begin{equation}
Y_1=q^{-1/2}t'_{11}t_{13}+t'_{12}t_{12}+q^{1/2}t'_{13}t_{11}
\end{equation}
and so on. The constraints involving $K''$ only through $X_i$ are
the following ones
\begin{eqnarray}
&&q^{1/2}\left(t_{11}t'_{31}-t'_{11}t_{31}\right)=\left(t_{21}t'_{21}-t'_{21}t_{21}\right)=
q^{-1/2}\left(t_{31}t'_{11}-t'_{31}t_{11}\right)=-K''X_{1},
\nonumber\\
&&q^{1/2}\left(t_{11}t'_{32}-t'_{11}t_{32}\right)=\left(t_{21}t'_{22}-t'_{21}t_{22}\right)=
q^{-1/2}\left(t_{31}t'_{12}-t'_{31}t_{12}\right)=-K''X_{2},
\nonumber\\
&&q^{1/2}\left(t_{12}t'_{31}-t'_{12}t_{31}\right)=\left(t_{22}t'_{21}-t'_{22}t_{21}\right)=
q^{-1/2}\left(t_{32}t'_{11}-t'_{32}t_{11}\right)=-K''X_{4},
\nonumber\\
&&q^{1/2}\left(t_{12}t'_{33}-t'_{12}t_{33}\right)=\left(t_{22}t'_{23}-t'_{22}t_{23}\right)=
q^{-1/2}\left(t_{32}t'_{13}-t'_{32}t_{13}\right)=-K''X_{6},
\nonumber\\
&&q^{1/2}\left(t_{13}t'_{32}-t'_{13}t_{32}\right)=\left(t_{23}t'_{22}-t'_{23}t_{22}\right)=
q^{-1/2}\left(t_{33}t'_{12}-t'_{33}t_{12}\right)=-K''X_{8},
\nonumber\\
&&q^{1/2}\left(t_{13}t'_{33}-t'_{13}t_{33}\right)=\left(t_{23}t'_{23}-t'_{23}t_{23}\right)=
q^{-1/2}\left(t_{33}t'_{13}-t'_{33}t_{13}\right)=-K''X_{9}.
\end{eqnarray}

There are six corresponding sets involving $K''$ only through
$Y_i$. As for (B.11) they are obtained by transposing indices in
the first three terms of each equation of (B.12) and changing the
sign before $K''$. Thus
\begin{equation}
q^{1/2}\left(t_{11}t'_{13}-t'_{11}t_{13}\right)=\left(t_{12}t'_{12}-t'_{12}t_{12}\right)
=q^{-1/2}\left(t_{13}t'_{11}-t'_{13}t_{11}\right)=K''Y_1
\end{equation}
and so on.

Finally there is a set involving $K''$ through both $X_i$ and
$Y_i$
\begin{eqnarray}
&&\left(t_{11}t'_{33}-t'_{11}t_{33}\right)=-K''q^{-1/2}\left(X_{3}-Y'_{3}\right),
\nonumber\\
&&\left(t_{12}t'_{32}-t'_{12}t_{32}\right)=-K''\left(q^{-1/2}X_{5}-Y'_{3}\right),
\nonumber\\
&&\left(t_{13}t'_{31}-t'_{13}t_{31}\right)=-K''\left(q^{-1/2}X_{7}-q^{1/2}Y'_{3}\right),
\nonumber\\
&&\left(t_{21}t'_{23}-t'_{21}t_{23}\right)=-K''\left(X_{3}-q^{-1/2}Y'_{5}\right),
\nonumber\\
&&\left(t_{22}t'_{22}-t'_{22}t_{22}\right)=-K''\left(X_{5}-Y'_{5}\right),
\nonumber\\
&&\left(t_{23}t'_{21}-t'_{23}t_{21}\right)=-K''\left(X_{7}-q^{1/2}Y'_{5}\right),
\nonumber\\
&&\left(t_{31}t'_{13}-t'_{31}t_{13}\right)=-K''\left(q^{1/2}X_{3}-q^{-1/2}Y'_{7}\right),
\nonumber\\
&&\left(t_{32}t'_{12}-t'_{32}t_{12}\right)=-K''\left(q^{1/2}X_{5}-Y'_{7}\right),
\nonumber\\
&&\left(t_{33}t'_{11}-t'_{33}t_{11}\right)=-K''q^{1/2}\left(X_{7}-Y'_{7}\right).
\end{eqnarray}
An alternative approach to the $\hat{R}tt$  relations is via the
diagonalization of $P_0'$. The diagonalizer is given in
\cite{R13}. Such an approach was presented for our multiparameter
("nested-sequence") class in App. C of Ref. 3.
\end{appendix}


\newpage

\begin{thebibliography}{99}
\itemsep=-.3pc

\bibitem{R1} A. Chakrabarti, {\it A nested sequence of
     projectors and corresponding braid matrices $\hat R(\theta)$: (1)
     odd dimensions}, Jour. Math. Phys. {\bf 46}, 063508 (2005)
\texttt{math.QA/0401207}.

\bibitem{R2} A. Chakrabarti, {\it Aspects of a new class of braid matrices: Roots of unity and hyperelliptic $q$
for triangularity, $L$-algebra, link-invariants, non-commutatives
spaces}, Jour. Math. Phys. {\bf 46}, 063509 (2005).

\bibitem{R3} B. Abdesselam and A. Chakrabarti, {\it A nested sequence of projectors: (2) Multiparameter
multistate statistical models, Hamiltonians, $S$-matrices}, Jour.
Math. Phys. {\bf 47}, 053508 (2006).

\bibitem{R4} R.J. Baxter, {\it Exactly solved models in statistical
mechanics}, Acad. Press (1982).

\bibitem{R5} P.P. Kulish and E.K. Sklyanin, {\it Integrable Quantum Field theories}, Lecture Notes in Physics (Springer,
New York, 1982) p. 61.

\bibitem{R6} H.J. De Vega, {\it Yang-Baxter algebras, integrable
theories and quantum groups}, Int. Jour. Mod. Phys. A vol. 4, 2371
(1989).

\bibitem{R7} L.D. Faddeev, {\it How algebraic Bethe ansatz works for
integrable models}, 1996 Les houches summer school 1995
[hep-th/9605187].

\bibitem{R8} D. Arnaudon, N. Cramp\'e, A. Doikou, L. Frappat and E.
Ragoucy, {\it Analytical Bethe Ansatz for closed and open
$gl(N)$-chains in any representation}, Jour. Stat. Mech., (2005)
02007.

\bibitem{R9} P. Isaev, Sov. J. Part. Nucl. 26, 501 (1995)

\bibitem{R10} J. Donin, P.P.Kulish and A.I. Mudrov, {\it On universal solution to reflection equation},
Lett. Math. Phys. 63, 179 (2003).

\bibitem{R11} V.E. Korepin, G. Izergin and N.M. Bogoliubov {\it
Quantum inverse scattering method, correlation functions and
algebraic Bethe Ansatz}, Cam. Univ. Press (1993).

\bibitem{R12} N. Kitanine, J.M. Maillet, N.A. Stanov and V. Terras,
{\it Spin spin correlation functions of the XXZ-1/2 Heisenberg
chain in a magnetic field}, Nucl. Phys. B 641, 487 (2002). {\it
Correlation functions of the XXZ spin-1/2 Heisenberg chain of the
free fermion point from their multiple integral representations},
Nucl. Phys. B 642, 433 (2002).

\bibitem{R13} A. Chakrabarti, {\it Canonical factorization and diagonalization of Baxterized braid matrices:
Explicit constructions and applications}, Jour. Math. Phys. {\bf
44}, 5320 (2003).

\end{thebibliography}
\end{document}